\documentclass[notitlepage,11pt,a4paper]{scrartcl}
\usepackage{amstext}
\usepackage{amssymb}
\usepackage{amsfonts}
\usepackage{amssymb}
\usepackage{amsmath}
\usepackage{graphicx}
\usepackage{mathrsfs}
\usepackage[arrow,matrix,curve]{xy}
\usepackage[latin1]{inputenc}
\usepackage{theorem}
\usepackage{geometry}
\newtheorem{satz}{Theorem}[section]
\newtheorem{defi}[satz]{Definition}
\newtheorem{bem}[satz]{Remark}
\newtheorem{kor}[satz]{Corollary}
\newtheorem{lem}[satz]{Lemma}
\newtheorem{ver}[satz]{Conjecture}

\newtheorem{bei}[satz]{Example}
\newtheorem{pro}[satz]{Proposition}

\newcommand{\modu}{~\mathrm{mod}~}

\newcommand{\qed}{$\hfill\square$\\}
\newcommand{\filt}[6]{
\[
\begin{xy}
\xymatrix@R20pt@C20pt{
&\mathbb{C}^3&\\\langle #1,#2 \rangle\ar[ru]&\langle #3,#4\rangle\ar[u]&\langle #5,#6\rangle\ar[lu]\\\langle #1\rangle\ar[u]&\langle #3\rangle\ar[u]&\langle #5\rangle\ar[u]}
\end{xy}
\]
} 
\begin{document}
  \title{Localization in quiver moduli spaces}
  \author{Thorsten Weist\\Fachbereich C - Mathematik\\
Bergische Universität Wuppertal\\
D - 42097 Wuppertal, Germany\\
\texttt{e-mail: weist@math.uni-wuppertal.de}}
\date{} 

  \maketitle
\begin{abstract}
Torus fixed points of quiver moduli spaces are given by stable representations of the universal (abelian) covering quiver. As far as the Kronecker quiver is concerned they can be described by stable representations of certain bipartite quivers coming along with a stable colouring. By use of the glueing method it is possible to construct a huge class of such quivers implying a lower bound for the Euler characteristic. For certain roots it is even possible to construct all torus fixed points. 

\end{abstract}
\maketitle
\section{Introduction}
\noindent A common method providing topological information of algebraic varieties is the consideration of fixed points under a torus action. For instance the Euler characteristic is already given by the Euler characteristic of their fixed point components. If we consider moduli spaces of stable quiver representations, we also often obtain interesting objects as fixed point components like indecomposable tree modules in the case of the Kronecker quiver, see \cite{wei3}. In general, torus fixed points of quiver moduli spaces are given by representations of the universal (abelian) covering quiver.\newline

The main focus of this paper is on torus fixed points of Kronecker moduli spaces. It is inasmuch particularly interesting as by use of the localization method we are able to prove parts of a conjecture based on Michael Douglas \cite{dou} concerning the Euler characteristic of these moduli spaces. It says that for coprime dimension vectors $(d,e)$ the logarithm of the Euler characteristic continuously depends on the fraction $\frac{e}{d}$. More specified this means that there exists a continuous function $f$ such that for every coprime dimension vector $(d,e)$ there exists another dimension vector $(d_s,e_s)$ such that
\[f\left(\frac{e}{d}\right)=\lim_{n\rightarrow\infty}\frac{\ln{\chi(M_{d_s+nd,e_s+ne}^m)}}{d_s+nd}.\]
In particular, the right hand side converges. In \cite{wei} a candidate for this function could be determined and it could be proved that under the assumption of continuity the function is already uniquely determined by $f(1)$.

Even if continuity is still an open question, by use of the localization method we are able to calculate the value at the point one. Actually, we are able to determine a formula for the Euler characteristic of the Kronecker moduli spaces for the dimension vectors $(d,d+1)$. Moreover, we show that the Euler characteristic grows exponentially with the dimension vector which is an immediate consequence of the conjecture.

The paper is organised as follows: in the second section the notion of quivers and their representations is introduced. Moreover, we recall the definition of stability and general results concerning the representation spaces which are needed in the sequel.\\

In the third section we consider a torus action on quiver moduli spaces. We show that the fixed points of moduli spaces of quivers without oriented cycles are exactly the stable representations of the universal abelian covering quiver. Actually, after suitable many localization steps the remaining torus fixed points are representations of the universal covering quiver.\\

In the fourth section we apply the localization method to the generalized Kronecker quiver. The universal covering quiver is a regular $m$-tree coming along with a certain orientation. In particular, every stable representation of a bipartite quiver which can be embedded into this $m$-tree is a torus fixed point. Therefore, we investigate stable bipartite quivers in more detail, i.e. quivers with a fixed dimension vector allowing at least one stable representation. We construct stable bipartite quivers of dimension type $(d_s,e_s)+n(d,e)$ by glueing certain bipartite quivers of dimension types $(d_s,e_s)$ and $(d,e)$. Thereby the dimension vector $(d_s,e_s)$ is uniquely determined by $(d,e)$. The dimension type of a bipartite quiver is given by the sum of the dimensions of the sources and sinks respectively. In this way, for every coprime dimension type we can construct a huge class of such quivers.\\

In the fifth section we briefly treat combinatorics of trees. With the stated methods we can asymptotically count the number of stable bipartite quivers constructed in the preceding section. In some cases we can even count it exactly.\\

In the last section several applications of the developed methods are treated. After investigating Douglas' conjecture in more detail, we study the mentioned function $f$ at the point one. Since all localization data of dimension type $(d,(m-1)d+1)$ are known, we can determine a formula for the Euler characteristic in this case. By applying the reflection functor we can also determine $f(1)$.

Afterwards, by use of the methods of the fourth section and combinatorics of trees we can determine a lower bound for the Euler characteristic for every coprime dimension vector. In particular, we prove that the Euler characteristic grows exponentially with the dimension vector.

In the fourth subsection the case of the dimension vector $(3,4)$ is given as an detailed example. The fifth subsection deals with the dimension vector $(d,d)$. We prove that there does not exist any stable representation of the universal covering quiver if $d\geq 2$ because torus fixed points of this dimension type are always cyclic. Thus it follows that the Euler characteristic vanishes in this case. In the sixth subsection we answer a question posed in \cite{dre}: when does there exist fixed point components containing infinitely many fixed points? Actually, there exist only finitely many torus fixed points for dimension vectors $(d,e)$ such that $d=1,2$ or in the associated reflected cases.
\section{Generalities}\label{allg}
\noindent Let $k$ be an algebraically closed field.
\begin{defi}
A quiver $Q=(Q_0,Q_1,t,h)$ is a quadruple consisting of a set of vertices $Q_0$, a set of arrows $Q_1$ and two maps $h,\,t:Q_1\rightarrow Q_0$ which associate its tail $t(\alpha)$ and its head $h(\alpha)$ to an arrow $\alpha\in Q_1$.

A vertex $q\in Q_0$ is called a sink (resp. a source) if $t^{-1}(q)=\emptyset$ (resp. $h^{-1}(q)=\emptyset$). 

A quiver is bipartite if $Q_0=I\cup J$ such that every vertex $i\in I$ is a source and every vertex $j\in J$ is a sink.
\end{defi}

We will also denote an arrow by $\alpha:q\rightarrow q'$ which means that $t(\alpha)=q$ and $h(\alpha)=q'$ for $q,\, q'\in Q_0$. In the following we only consider quivers without oriented cycles. 

Define the abelian group
\[\mathbb{Z}Q_0=\bigoplus_{q\in Q_0}\mathbb{Z}q\] and the monoid of dimension vectors $\mathbb{N}Q_0\subset \mathbb{Z} Q_0$. On $\mathbb{Z}Q_0$ we define a (non-symmetric) bilinear form, called the Euler form, by 
\[\langle d,e\rangle:=\sum_{q\in Q_0}d_qe_q-\sum_{\alpha\in Q_1}d_{t(\alpha)}e_{h(\alpha)}\]
for $d,\,e\in\mathbb{Z}Q_0$. Moreover, let $\{d,e\}:=\langle d,e\rangle+\langle e,d\rangle$ be the symmetrized Euler form.

A finite-dimensional $k$-representation of $Q$ is given by a tuple
\[X=((X_q)_{q\in Q_0},(X_{\alpha})_{\alpha\in Q_1}:X_{t(\alpha)}\rightarrow X_{h(\alpha)})\]
of finite-dimensional $k$-vector spaces and $k$-linear maps between them. The dimension vector $\underline{\dim}X\in\mathbb{N}Q_0$ of $X$ is defined by
\[\underline{\dim}X=\sum_{q\in Q_0}\dim_kX_qq.\]
The support $\mathrm{supp}(d)$ of a dimension vector $d\in\mathbb{N}Q_0$ is the quiver defined by the vertices $\mathrm{supp}(d)_0=\{q\in Q_0\mid d_q\neq 0\}$ and the arrows $\mathrm{supp}(d)_1=\{\alpha\in Q_1\mid h(\alpha),\,t(\alpha)\in \mathrm{supp}(d)_0\}$. In the following, we only consider dimension vectors which support is finite, i.e. the number of vertices and arrows is finite.

Let $d\in\mathbb{N}Q_0$ be a dimension vector. The variety $R_d(Q)$ of $k$-representations of $Q$ with dimension vector
$d$ is defined as the affine $k$-space
\[R_d(Q)=\bigoplus_{\alpha\in Q_1} \mathrm{Hom}_k(k^{d_{t(\alpha)}},k^{d_{h(\alpha)}}).\]
The algebraic group $G_d=\prod_{q\in Q_0} Gl_{d_q}(k)$
acts on $R_d(Q)$ via simultaneous base change, i.e.
\[(g_q)_{q\in Q_0}\ast
(X_{\alpha})_{\alpha\in Q_1}=(g_{h(\alpha)}X_{\alpha}g_{t(\alpha)}^{-1})_{\alpha\in Q_1}.\] The orbits are in bijection with the isomorphism classes of
$k$-representations of $Q$ with dimension vector $d$.\\

In the space of $\mathbb{Z}$-linear functions $\mathrm{Hom}_{\mathbb{Z}}(\mathbb{Z}Q_0,\mathbb{Z})$ we consider the basis given by the elements $q^{\ast}$ for $q\in Q_0$, i.e.
$q^{\ast}(q')=\delta_{q,q'}$ for $q'\in Q_0$. Define $\dim:=\sum_{q\in Q_0}q^{\ast}.$
After choosing $\Theta\in
\mathrm{Hom}_{\mathbb{Z}}(\mathbb{Z}Q_0,\mathbb{Z})$,
we define the slope function
$\mu:\mathbb{N}Q_0\backslash\{0\}\rightarrow\mathbb{Q}$ via
\[\mu(d)=\frac{\Theta{(d)}}{\dim(d)}.\]
The slope $\mu(\underline{\dim}X)$ of a representation $X\neq 0$ of $Q$ is abbreviated to $\mu(X)$.
\begin{defi}
A representation $X$ of $Q$ is semistable (resp. stable) if for all proper subrepresentations  $0\neq U\subsetneq X$ the following holds:
\[\mu(U)\leq\mu(X)\text{ (resp. } \mu(U)<\mu(X)).\]
\end{defi}

This definition is equivalent to that of A.King \cite{kin}: let $\tilde{\Theta}\in\mathrm{Hom}(\mathbb{Z}Q_0,\mathbb{Z})$ be a linear form. A representation $X$ is $\tilde{\Theta}$-semistable (resp. $\tilde{\Theta}$-stable) in the sense of King if $\tilde{\Theta}(\underline{\dim}X)=0$ and
\[\tilde{\Theta}(\underline{\dim} U)\geq 0\,\,(\mathrm{resp.}\,\, \tilde{\Theta}(\underline{\dim} U)>0)\] 
for all proper subrepresentations $0\neq U\subsetneq X$. Now fixing a representation $X$ define 
\[\tilde{\Theta}:=\mu(X)\cdot \dim-\Theta.\]
It is easy to check that the representation $X$ is semistable (resp. stable) in the former sense if and only if it is $\tilde{\Theta}$-semistable (resp. $\tilde{\Theta}$-stable).

Denote the set of semistable (resp. stable) points by $R^{ss}_d(Q)$ (resp. $R^s_d(Q)$). In this situation we have the following theorem going back to Mumford's GIT and which was proved by King, see \cite{mum}, \cite{kin}:
\begin{satz}\label{kingsatz}
\begin{enumerate}
\item The set of stable points $R^s_d(Q)$ is an open subset of the set of semistable points $R^{ss}_d(Q)$, which is an open subset of $R_d(Q)$. \item There exists a categorical quotient $M^{ss}_d(Q):=R^{ss}_d(Q)//G_d$. Moreover, $M^{ss}_d(Q)$ is a projective variety. \item There exists a geometric quotient $M^{s}_d(Q):=R_d^{s}(Q)/G_d$, which is a smooth open subvariety of $M^{ss}_d(Q)$.
\end{enumerate}
\end{satz}
Note that the set of semistable (resp. stable) points of $R_d(Q)$ can be empty. For a detailed description of the theory of quotients see \cite{muk}. We just briefly treat the construction. All statements about algebraic groups applied in this paper can for instance be found in \cite{spr} or \cite{hum}.

We obtain the quotient $M^{ss}_d(Q)$ called moduli space in what follows by defining a character $\chi$ of $G_d$ by
\[\chi((g_q)_{q\in Q_0}):=\prod_{q\in
Q_0}\det(g_q)^{\Theta(d)-\dim d\cdot\Theta_q},\] where $\Theta$ is the linear form obtained from the previous consideration.

For an affine variety $X$ endowed with an action of a reductive algebraic $G$ group the set of semi-invariants of weight $\chi^n$ is defined by
\begin{center}
$k[X]^{G,\chi^n}:=\{f\in k[X]\mid f(g\ast x)=\chi(g)^n\cdot
f(x)\;\forall\, g\in G,\,\forall\, x\in X\}.$
\end{center}
Furthermore, the ring of $\chi$-semi-invariants is given by
\begin{center}
$k[X]^G_{\chi}:=\bigoplus\limits_{n=0}^\infty k[X]^{G,\chi^n}$.
\end{center}
Then we have
\[M^{ss}_d(Q)=\mathrm{Proj}(k[X]^G_{\chi}),\] the projective spectrum of the ring of semi-invariants.
\begin{bem}\label{bem1}
\end{bem}
\begin{enumerate}
\item Since there exists only one closed orbit, the affine quotient is just a point. Therefore, we get $k[R_d(Q)]^G=k$. Thus the projective quotient has no affine component and is a projective variety.
\item Since $R_d(Q)$ is an affine space and thus smooth, we get that the open subset of stable points is smooth. Thus, since the moduli space
$M^s_d(Q)$ is an orbit space associated to the group action restricted to the stable points, it is smooth as well.
\item The moduli space $M_d^{ss}(Q)$ does not parametrize the semistable representations, but the polystable ones. Polystable representations are such representations which can be decomposed into a direct sum of stable ones of the same slope.
\item If semistability and stability coincide, $M_d^{ss}(Q)$ actually is a smooth projective variety. Obviously this is the case if $\mu(d)\neq\mu(e)$ for all $0\neq e<d$. In this case the dimension vector $d$ is said to be $\Theta$-indivisible.
\item
For a stable representation $X$ we have that its orbit is of maximal possible dimension. Since the scalar matrices act trivially on $R_d(Q)$, the isotropy group is at least one-dimensional. Thus, if the moduli space $M^s_d(Q)$ is not empty, for the dimension of the moduli space we have \[\dim M^{s}_d(Q)=1-\langle d, d\rangle.\]
\end{enumerate}

Finally we point out some properties of (semi-)stable representations. These properties will be very useful at different points of this paper, for proofs see \cite{hn}. \begin{lem}\label{zusa} For a quiver $Q$ let $0\rightarrow
M\rightarrow X\rightarrow N\rightarrow 0$ be a short exact sequence of representations.\begin{enumerate}
\item The following are equivalent:
\begin{enumerate}
\item $\mu(M)\leq\mu(X)$
\item $\mu(X)\leq\mu(N)$
\item $\mu(M)\leq\mu(N)$
\end{enumerate}
\item The following holds: $\min(\mu(M),\mu(N))\leq\mu(X)\leq \max(\mu(M),\mu(N)).$
\item If $\mu(M)=\mu(X)=\mu(N)$, then $X$ is semistable if and only if $M$ and $N$ are semistable.
\end{enumerate}
\end{lem}

From the first property we immediately get that stable representations are indecomposable.
Denote by $E_q$ the simple representation corresponding to the vertex $q$ defined by $(E_q)_{q}=k$ and $(E_q)_{q'}=0$ for $q'\in Q_0$ with $q'\neq q$.

Define $e_q\in\mathbb{Z}Q_0$ by $(e_q)_{q'}:=\delta_{q,q'}$ which is the dimension vector of $E_q$.
For a quiver $Q$ consider the matrix $A=(a_{q,q'})_{q,q'\in Q_0}$ defined by $a_{q,q'}=\{e_q,e_{q'}\}$ for $q,\,q'\in Q_0$. Fixing some $q\in Q_0$ define $r_q:\mathbb{Z}Q_0\rightarrow\mathbb{Z}Q_0$ by
\[r_q(e_{q'})=e_{q'}-a_{q,q'}\cdot e_q.\]
Let $Q_q$ be the quiver resulting from $Q$ by reversing all arrows with head or tail $q$. We have the following theorem, see \cite{bgp}:
\begin{satz}\label{kspi}
Let $Q$ be a quiver and $q\in Q_0$ a fixed vertex. Let $q$ be a sink (resp. a source). Then there exists a functor \[R_q^+ (\text{resp. } R_q^-):\modu kQ\rightarrow \modu kQ_q\] with the following properties (if $q$ is a source, replace $+$ by $-$):
\begin{enumerate}
\item $R_q^+(U\oplus U')=R_q^+(U)\oplus R_q^+(U')$.
\item Let $U$ be an indecomposable representation of $Q$.
\begin{enumerate}
\item If $U\cong E_q$, then $R_q^+(E_q)=0$.
\item If $U\ncong E_q$, then $R_q^+(U)$ is indecomposable with $R_q^-R_q^+(U)\cong U$ and we have $\underline{\dim} R_q^+(U)=r_q(\underline{\dim}(U))$.
\end{enumerate}
\end{enumerate}
Moreover, we have: $\mathrm{End} U\cong \mathrm{End} R_q^+(U)$.
\end{satz}
\section{Localization in quiver moduli spaces}
\noindent Analogously to \cite{rei3}, in this section we introduce the localization in moduli spaces of stable representations. Some of the ideas are based on localization techniques in moduli spaces of simple representation provided by \cite{rei}. An explicit method to detect fixed points of these moduli spaces under a torus action is explained. These fixed points are stable representations of the universal abelian covering quiver.
\subsection{Torus fixed points}\label{univcov}
\noindent For the remaining part of the paper we fix $k=\mathbb{C}$. Let $G$ be an algebraic group and $\chi:G\rightarrow\mathbb{C^{\ast}}$ be a character of $G$, i.e. a morphism of algebraic groups. Denote by $X(G)$ the set of all characters of $G$ with the group structure given in the obvious way. In the following the composition is written additively.

Let $G$ be a closed subgroup of $Gl(V)$ and $V$ be a representation of $G$. For all characters $\chi\in X(G)$
define the semi-invariants of weight $\chi$ by
\[V_{\chi}=\{v\in V\mid g\cdot v=\chi(g)v \;\forall g\in G\}.\]
Note that if $\varphi:G\rightarrow Gl(V)$ is a rational representation, the definition can be transferred. It is well known that we obtain a decomposition into weight spaces $V=\bigoplus_{\chi\in X(G)}V_{\chi}$.\\

Let further $T:=(\mathbb{C}^{\ast})^{|Q_1|}$ be the $|Q_1|$-dimensional torus. It acts on $R_d(Q)$ via
\[((t_{\alpha})_{\alpha\in Q_1})\cdot((X_{\alpha})_{\alpha\in
Q_1})=(t_{\alpha}\cdot X_{\alpha})_{\alpha\in Q_1}.\] 
Since the torus action commutes with the $G_d$-action, it induces a $T$-action on $M^s_d(Q)$.

Since the scalar matrices act trivially on $R_d(Q)$, the
$G_d$-action factorises through the quotient $PG_d:=G_d/\mathbb{C}^{\ast}.$
Let $X\in M^s_d(Q)$ be a fixed point under the torus action. Considering the algebraic group
\[G:=\{((g_q)_{q\in Q_0},t)\in PG_d\times T\mid t\cdot X=(g_q)_{q\in Q_0}\ast X\}\]
we get projections $p_1:G\rightarrow PG_d$ and $p_2:G\rightarrow T$ respectively with the following property:
\begin{lem}
Let $X$ be a torus fixed point. The following holds:
\begin{enumerate}
\item The projection $p_2:G\rightarrow T$ is an isomorphism. \item
In particular, the projection $p_1:G\rightarrow PG_d$ induces a homomorphism of algebraic groups $\varphi:=p_1\circ
p_2^{-1}:T\rightarrow PG_d$ such that $\varphi(t)\ast X=t\cdot X.$
\end{enumerate}
\end{lem}
$\it{Proof.}$ Since $X\in M^s_d(Q)$ is a fixed point, $p_2$
is surjective. Moreover, since $X$ is stable, its orbit is of maximal possible dimension. Thus the isotropy group of $X$ under the action of
$PG_d$ is trivial implying the injectivity. The second part immediately follows from this.\qed

By the preceding considerations we get:
\begin{lem}\label{fixp}
The following are equivalent:
\begin{enumerate}
\item $X$ is a fixed point.
\item There exists a morphism of algebraic groups $\varphi:T\rightarrow
PG_d$ such that
\[(t_{\alpha})_{\alpha\in Q_1}\cdot (X_{\alpha})_{\alpha\in Q_1}=\varphi((t_{\alpha})_{\alpha\in
Q_1})\ast (X_{\alpha})_{\alpha\in Q_1}\] for all
$(t_{\alpha})_{\alpha\in Q_1}\in T$.
\end{enumerate}
\end{lem}

The following lemma assures that we get a weight space decomposition of the vector space corresponding to some fixed point:
\begin{lem}\label{lift}
Let $T\cong(\mathbb{C}^{\ast})^m$ with $m\geq 1$ be a torus. Every homomorphism of algebraic groups $\varphi:T\rightarrow PGl_d(\mathbb{C})$ can be lifted, i.e. there exists a homomorphism of algebraic groups $\psi:T\rightarrow Gl_d(\mathbb{C})$ such that $\varphi=\pi\circ\psi$.
\end{lem}
{\it Proof.}
In general, if $G$ is a reductive algebraic group, the image of a morphism $\varphi:T\rightarrow G$ is again a torus and, therefore, contained in a maximal torus $T_0\subset G$. Since all maximal tori are conjugate, we can assume that $T_0=(\mathbb{C}^{\ast})^n$ for some $n\in\mathbb{N}$. Thus, in order to prove the statement, it suffices to prove that every morphism $\varphi:\mathbb{C}^{\ast}\rightarrow(\mathbb{C}^{\ast})^{n-1}$ can be lifted to a morphism $\psi:\mathbb{C}^{\ast}\rightarrow(\mathbb{C}^{\ast})^n$ where $\pi:(\mathbb{C}^{\ast})^n\rightarrow(\mathbb{C}^{\ast})^{n-1}$ is the projection induced by the projection $\pi:Gl_d(\mathbb{C})\rightarrow PGl_d(\mathbb{C})$. Note that, since $\pi$ is surjective, every maximal torus is mapped to a maximal torus. Now if $\varphi(t)=(t^{r_2},\ldots,t^{r_n})$ and $\pi(t_1,\ldots,t_n)=(\frac{t_2}{t_1},\ldots,\frac{t_n}{t_1})$, then we may set $\psi(t)=(1,t^{r_2},\ldots,t^{r_n})$.\qed

The lift $\psi:T\rightarrow G_d$ for $\varphi$ can be decomposed in $|Q_0|$ morphisms of algebraic groups $\psi_q:T\rightarrow Gl_{d_q}$. Thus if $X$ is a fixed point and $\varphi:T\rightarrow PG_d$ the corresponding morphism, we can fix a lift in order to get a weight space decomposition of each vector space, i.e. $X_q=\bigoplus_{\chi\in X(T)}X_{q,\chi}$. For a $d$-dimensional torus we have $X(T)\simeq\mathbb{Z}^d$. As far as the torus $(\mathbb{C}^{\ast})^{|Q_1|}$ is concerned we denote the canonical basis of $X(T)$ by $(e_{\alpha})_{\alpha\in Q_1}$.  The next lemma shows that these weight space decompositions are compatible with the linear maps:
\begin{lem}
Let $X=(X_{\alpha})_{\alpha\in Q_1}$ be a fixed point under the torus action. Let $\varphi:T\rightarrow PG_d$ be the corresponding morphism of algebraic groups and fix a lift $\psi:T\rightarrow G_d$. Then we have:
\[X_{\alpha}(X_{q,\chi})\subseteq X_{q',\chi+e_{\alpha}}\text{ for all
}\chi\in X(T),\,\alpha:q\rightarrow q'.\]

\end{lem}
$\it{Proof.}$ Let $t=(t_{\alpha})_{\alpha\in Q_1}\in T$ and
$x\in X_{q,\chi}$. Then we have
$$\psi_{q'}(t)X_{\alpha}(x)=\psi_{q'}(t)
X_{\alpha}\psi_q(t)^{-1}\psi_q(t)(x)=t_{\alpha}
X_{\alpha}\chi(t)(x)=(\chi+e_{\alpha})(t)X_{\alpha}(x).$$
\qed\\

Now we investigate the stability criterion for fixed points. We will see that it is enough to consider subspaces compatible with the weight space decomposition. This is important for the practicability of the introduced construction.

Let $X$ be a quiver representation. Define by $\mathrm{scss}(X)$
(\textbf{s}trongly \textbf{c}ontra\-dicting
\textbf{s}emi\-\textbf{s}tability) the subrepresentation $Y\subset
X$ for which the following holds:
\begin{enumerate}
\item $\mu(Y)=\max\{\mu(U)\mid U\subset X\}.$
\item $\dim(Y)=\max\{\dim(U)\mid U\subset X, \mu(U)=\mu(Y)\}.$ 
\end{enumerate}
Thus $Y$ is of maximal dimension among the subrepresentations with maximal slope. It is straightforward to check that the subrepresentation $\mathrm{scss}(X)$ is uniquely determined, see e.g. \cite{rei3}.

\begin{lem}\label{kompa}
Let $X$ be a fixed point under the torus action with dimension vector $d$.
Let \[X_q=\bigoplus_{\chi\in X(T)}X_{q,\chi}\] be the weight space decomposition with respect to some fixed lift $\psi$ of the associated morphism
$\varphi:T\rightarrow PG_d$. Then the following are equivalent:
\begin{enumerate}
\item $X$ is semistable (resp. stable).
\item For all subrepresentations $U$, which are compatible with the weight space decomposition of $X$, i.e.
$U_q=\bigoplus_{\chi\in X(T)}U_{q,\chi}$ for all $q\in Q_0$ where
$U_{q,\chi}\subset X_{q,\chi}$, we have $\mu(U)<\mu(X)$ (resp. $\mu(U)\leq\mu(X)$).
\end{enumerate}
\end{lem}
{\it Proof.} One conclusion is clear. Thus let $X$ be a representation satisfying the second property.  We first show that $X$ is semistable.

Let $U=\mathrm{scss}(X)$. Since $X$ is a fixed point, by Lemma \ref{fixp} there exists a morphism of algebraic groups $\varphi:T\rightarrow PG_d$ such that \[(t_{\alpha})_{\alpha\in Q_1}\cdot (X_{\alpha})_{\alpha\in Q_1}=\varphi((t_{\alpha})_{\alpha\in
Q_1})\ast (X_{\alpha})_{\alpha\in Q_1}\] for all
$(t_{\alpha})_{\alpha\in Q_1}\in T$. By Lemma \ref{lift} we may choose a lift $\psi:T\rightarrow G_d$. Fix such a lift and consider
\[\psi(t)U:=(\psi_q((t_{\alpha})_{\alpha\in
Q_1})(U_q))_{q\in Q_0}\]for each $(t_{\alpha})_{\alpha\in
Q_1}\in T$. 

Hence for each arrow $\alpha:q\rightarrow q'$ we obtain \[X_{\alpha}\psi_q(t)U_q=\frac{1}{t_{\alpha}}\psi_{q'}(t)X_{\alpha}\psi_q(t)^{-1}\psi_q(t)U_q\subset
\psi_{q'}(t)U_{q'}\] because $X_{\alpha}U_q\subset U_{q'}$. Thus
$\psi(t)U$ is a subrepresentation of $X$. Since $\psi_q(t)$
is invertible for every $q\in Q_0$, the dimension vectors of $U$ and $\psi(t)U$ coincide. Because of the uniqueness of
$\mathrm{scss}(X)$ it follows that $\psi(t)U=U$ for all $t\in T$. This is equivalent to $\psi_q(t)U_q=U_q$ for all $t\in T$ and
all $q\in Q_0$. This implies that $U=\mathrm{scss}(X)$ is 
compatible with the weight space decomposition. Therefore, by assumption we have that $\mu(\mathrm{scss}(X))\leq\mu(X)$. Hence $X$ is semistable. Indeed, the slope of $\mathrm{scss}(X)$ is maximal among the set of subrepresentations of $X$.\\

Thus it remains to show that if proper inequality holds, it follows that $X$ is stable. Assume that $X$ is not stable. By the preceding considerations we know that $X$ is semistable. Thus we may assume that there exists a subrepresentation $U$ of $X$ such that $\mu(U)=\mu(X)$. Consider again the lift $\psi:T\rightarrow G_d$ from above. As above we obtain that $\{\psi(t)U\mid t\in T\}$ is a set of subrepresentations of $X$. Let $e:=\underline{\dim} U$ and consider the quiver Grassmannian $\hat{\mathrm{Gr}}_e(X)$ of subrepresentations of dimension $e$ of $X$ which is a projective variety because it is a closed subvariety of the product of the usual Grassmannians $\mathrm{Gr}_{e_q}(X_q)$, $q\in Q_0$. It is also not empty because $U\in\hat{\mathrm{Gr}}_e(X)$. From the considerations above we obtain a torus action on $\hat{\mathrm{Gr}}_e(X)$ given by $(t,U)\mapsto\psi(t)\cdot U$. But since $\hat{\mathrm{Gr}}_e(X)$ is projective it follows by Borel's Fixed Point Theorem that the fixed point set is not empty. Thus there exists a subrepresentation $U'\subset X$ with $\mu(U')=\mu(X)$ such that $\psi(t)\cdot U'=U'$. But this again means that $U'$ is compatible with the weight space decomposition which is contradiction.
\qed

If we choose another lift $\psi '$, one easily verifies that there exists a character $\chi\in X(T)$ such that $\psi=\chi\psi'.$ If we have two representatives $X$ and $X'$ of a fixed point, i.e. there exists some $g\in G_d$ such that $X'=g\ast X,$ we can assume that the weight space decomposition does not change. Indeed, if $\varphi$ is the morphism belonging to $X$, for the morphism $\varphi'$ belonging to $X'$ we have $\varphi'=\pi(g)\cdot\varphi\cdot \pi(g^{-1})$ where $\pi:G_d\rightarrow PG_d$ is the canonical projection. Thus, if $\psi$ is a lift of $\varphi$, we have that $\psi'=g\cdot\psi\cdot g^{-1}$ is a lift of $\varphi'$. Thus we obtain that the dimensions of the weight spaces for both morphism $\psi$ and $\psi'$ coincide. Indeed, for $x\in X_{q,\chi}$ we have $\chi(t)x=\psi(t)_qx$ if and only if $\chi(t)g_qx=\psi'(t)_qg_qx$.  This also shows that $g$ is compatible with the weight space decomposition.\\

We want to define a quiver such that the fixed point components correspond to moduli spaces of this quiver with compatible dimension vectors. Therefore, define the quiver $\hat{Q}$ by the vertex set
\[\hat{Q}_0=Q_0\times X(T)\]
and for each arrow $\alpha:q\rightarrow q'$ and each character $\chi\in X(T)$ we have an arrow
\[(\alpha,\chi):(q,\chi)\rightarrow (q',\chi + e_{\alpha})\]
in $\hat{Q}_1$. This is the universal abelian covering quiver of
$Q$.

Let $X$ be a fixed point of the moduli space with respect to the torus action. Then define the corresponding dimension vector $\hat{d}\in\mathbb{N}\hat{Q}_0$
by \[\hat{d}_{q,\chi}:=\dim_{\mathbb{C}}X_{q,\chi}.\] Obviously $X$ can be considered as representation of this quiver.

The stability condition for representations of this quiver is induced from $\Theta$, i.e. we define a linear form $\hat{\Theta}:\mathbb{Z}\hat{Q}_0\rightarrow\mathbb{Z}$ such that
\[\hat{\Theta}_{q,\chi}=\Theta_q\]
for all $q\in Q_0$ and all $\chi\in X(T)$. Thus by Lemma $\ref{kompa}$, semistable (resp. stable) fixed points can be identified with semistable (resp. stable) representations of the just introduced quiver.\\

Next we show that such a representation corresponding to a fixed point $X$ is unique in a certain way. By the preceding considerations choosing another lift
$\psi$ just changes the weights of the weight space decomposition by translation by some character $\mu$. This corresponds to a group action of $\mathbb{Z}Q_1$ on $\hat{Q}_0$ defined by
\[\mu\cdot (q,\chi)=(q,\chi + \mu).\]
Now this induces a group action on the set of dimension vectors $\mathbb{N}\hat{Q}_0$. Two dimension vectors contained in the same orbit are said to be equivalent in the following. In the following, we consider the dimension vectors of $\hat{Q}$ up to this equivalence. Thus we have in conclusion:
\begin{satz}
For all fixed points $X\in M_d^s(Q)^T$ there exists (up to equivalence)
a unique dimension vector $\hat{d}$ for $\hat{Q}$ such that
$X$ corresponds to a stable representation of $\hat{Q}$ with dimension vector $\hat{d}$.
\end{satz}
\subsection{Description of fixed points}
\noindent Converse to the last section we construct an embedding of stable representations of the quiver $\hat{Q}$ into the fixed point set of the related moduli space. Therefore, fixing a representation of $\hat{Q}$ we construct a representation of $Q$ and show that the latter one is a fixed point.\\

Again consider the quiver $\hat{Q}$ and let $\hat{d}$ be a dimension vector. Define $d_q$ with $q\in Q_0$ by 
\[d_q=\sum_{\chi\in X(T)}\hat{d}_{q,\chi}.\]
We call a dimension vector $\hat{d}$ satisfying this property compatible with $d:=(d_q)_{q\in Q_0}$.

Let $\hat{X}=((\hat{X}_{q,\chi})_{q\in Q_0,\chi\in X(T)},(\hat{X}_{\alpha,\chi})_{\alpha\in Q_1,\chi\in X(T)})$ be a representation of $\hat{Q}$. Define a representation $X$ of $Q$ by the vector spaces
\[X_q:=\bigoplus_{\chi\in X(T)} \hat{X}_{q,\chi}.\]
and the linear maps 
\[X_{\alpha}=\bigoplus_{\chi\in X(T)} \left(\hat{X}_{\alpha,\chi}:\hat{X}_{q,\chi}\rightarrow\hat{X}_{q',\chi+e_{\alpha}}\right)\]
for all $\alpha:q\rightarrow q'$.

This defines a linear map $P:R_{\hat{d}}(\hat{Q})\rightarrow R_d(Q).$
Moreover, an embedding of $G_{\hat{d}}$ in $G_d$ arises from the decomposition of the vector spaces $X_q$ for $q\in Q_0$. Since the linear map is equivariant under the group action of $G_{\hat{d}}$, the map $P$ induces a map \[P:M^{ss}_{\hat{d}}(\hat{Q})\rightarrow M^{ss}_d(Q)\]
by use the universal property of the quotient. Furthermore, define a morphism of algebraic groups
$\psi=(\psi_q)_{q\in Q_0}:T\rightarrow G_d$ such that
$\psi_q:T\rightarrow Gl(X_q)$ is defined by
\[\psi_q(t)x=\chi(t)x\]
for all $t\in T$ and all $x\in X_{q,\chi}.$ This makes $\psi$ well-defined and by Lemma \ref{fixp} we obtain a morphism of algebraic groups $\varphi:T\rightarrow PG_d$ such that $P(X)=Y$ is a fixed point. Because of Lemma \ref{kompa} semistable (resp. stable) representations of $\hat{Q}$ are mapped to semistable (resp. stable) representations of $Q$. \begin{lem}
Let $X$ and $X'$ be stable representations of $\hat{Q}$ such that
$P(X)$ and $P(X')$ are isomorphic. Then $X$ and $X'$
are already isomorphic. \end{lem} {\it Proof.} Let $Y=P(X)$ and $Y'=P(X')$, define $d:=\underline{\dim}(Y)$
and let
\[g=(g_q\in Gl_{d_q}(\mathbb{C}))_{q\in Q_0}\]
be an isomorphism between $Y$ and $Y'$. We have
\[Y_{\alpha}'g_q=g_{q'}Y_{\alpha}\]
for all $\alpha:q\rightarrow q'\in Q_1$. We choose the morphism of algebraic groups $\psi=(\psi_q)_{q\in Q_0}$ corresponding to $Y$ as above. Since $Y'$ is a fixed point isomorphic to $Y$, by the considerations of the last subsection we can choose the lift $\psi'$ corresponding to $Y'$ such that we have $\psi_q'=g_q\psi_q g_q^{-1}$ for all $q\in Q_0$. But as before for these lifts we have
\[y\in Y_{q,\chi},\text{ i.e. } \psi(t)_qy=\chi(t)y\Leftrightarrow\,\psi'(t)_q(g_qy)=\chi(t)(g_qy).\]
Thus each $g_q$ induces an isomorphism between the weight spaces $Y_{q,\chi}$ and $Y'_{q,\chi}$. Hence we may understand $g$ as an isomorphism between $X$ and $X'$ because $g$ is compatible with the weight space decomposition.\qed



Every fixed point arises from such an embedding. Moreover, the images of these embeddings are pairwise disjoint so that we obtain the following concluding theorem:
\begin{satz}\label{isom}
The set of fixed points $M_d^s(Q)^T$ is isomorphic to the disjoint union of moduli spaces
\[\bigcup_{\hat{d}} M^s_{\hat{d}}(\hat{Q}),\] where $\hat{d}$ ranges over all equivalence classes of dimension vectors being compatible with $d$.
\end{satz}
\subsection{Euler characteristic of moduli spaces}
\noindent In this section we point out some basic properties of the Euler characteristic. For basics of algebraic topology see for instance \cite{lue}.

 Let $X$ be a smooth quasi-projective variety over the complex numbers of dimension $n$ and let $H^i(X)$, $i\in\mathbb{N}_0$, be the $i$-th singular cohomology group with coefficients in $\mathbb{C}$ which are $\mathbb{C}$-vector spaces satisfying $H^i(X)=0$ if $i>2n$ as is known. Define
\[h^i(X)=\dim_{\mathbb{C}}H^i(X).\]The Euler characteristic $\chi$ of $X$ is defined by
\[\chi(X)=\sum_{k=0}^{2n}(-1)^kh^k(X).\]
By the following theorem, which is a consequence of \cite[Chapter 2.5]{cg}, it follows that the localization method is suitable to calculate the Euler characteristic of varieties.\begin{satz}
Let $X$ be a smooth complex variety with a torus $T$ acting on it. Let $X^T$ be the fixed point set of $X$ under this action. Then for the Euler characteristic we have
\[\chi(X)=\chi(X^T).\]
\end{satz}

Note that this theorem also holds for complex varieties in general considering cohomology with compact support. Moreover, if the variety is smooth, the Euler characteristic in singular cohomology and the one in cohomology with compact support coincide.

By Theorem $\ref{isom}$ and because of the additivity of the Euler characteristic we obtain the following important result: \begin{satz}\label{zerlegung} Let $Q$ be a quiver with dimension vector $d$. Then for the Euler characteristic of the moduli space $M^s_d(Q)$ we have
\[\chi(M^s_d(Q))=\sum_{\hat{d}}\chi(M^s_{\hat{d}}(\hat{Q})),\]
where $\hat{Q}$ is the universal abelian covering quiver and $\hat{d}$
ranges over all equivalence classes being compatible with $d$.
\end{satz}

Let $Q$ be a quiver and $d$ be a $\Theta$-indivisible dimension vector. Consider the moduli space of stable representations $M^s_d(Q)$. From the formula for the Poincaré polynomials stated in \cite{rei2} we obtain that the coefficients corresponding to the monomials in odd degree vanish so that the odd cohomology vanishes. Moreover, from the Hard Lefschetz Theorem, see for instance \cite{gh}, we can conclude that
\[h^k(M^s_d(Q))\leq h^{k+2}(M^s_d(Q))\]
for $k<n$ and
\[h^k(M^s_d(Q))\geq h^{k+2}(M^s_d(Q))\]
for $k>n$ where $n$ is the dimension of $M^s_d(Q)$. Since we also have
\[h^0(M^s_d(Q))=h^{2n}(M^s_d(Q))=1,\]
we get the following result:
\begin{kor}
For moduli spaces of stable representations of a quiver $Q$ with $\Theta$-indivisible dimension vector $d$ we have:
\[\chi(M^s_d(Q))\geq \dim M^s_d(Q)+1.\]
\end{kor}
\subsection{Maps between universal quivers}
\noindent In this subsection we introduce the universal covering quiver of a connected quiver $Q$. Moreover, we construct maps from this quiver to the universal abelian covering quivers which are obtained by applying the localization technique recursively. Since these maps become injective on finite subquivers after finitely many steps, the remaining torus fixed points do not have a cyclic support.\\

Let $\mathfrak{Q}_1=\{\alpha,\alpha^{-1}\mid\alpha\in Q_1\}$ where $\alpha^{-1}$ is the formal inverse of $\alpha$. We will write $\alpha^{-1}:q'\rightarrow q$ for $\alpha:q\rightarrow q'\in Q_1$. A path $p$ is a sequence $(q_1\mid\alpha_1\alpha_2\ldots\alpha_n\mid q_{n+1})$ such that $\alpha_i:q_i\rightarrow q_{i+1}\in \mathfrak{Q}_1$. Thereby, we have the equivalence generated by
\[(q\mid\alpha\alpha^{-1}\mid q)\sim (q\mid\mid q).\]
In what follows, we always consider paths up to this equivalence. The set of words in $Q$ is generated by the arrows and their formal inverses, i.e. for a word $w$ we have $w=\alpha_1\ldots\alpha_n$ where $\alpha_i\in \mathfrak{Q}_1$. Denote the set of words of $Q$ by $W(Q)$. Note the difference between paths and words, i.e. a word may consist of any concatenation of arrows and their formal inverse whence two paths can only be concatenated if the head of one of the paths coincides with the tail of the other one.
The universal covering quiver $\tilde{Q}$ of $Q$ is given by the vertex set 
\[\tilde{Q}_0=\{(q,w)\mid q\in Q_0,w\in W(Q)\}\]
and the arrow set
\[\tilde{Q}_1=\{\alpha_{(q,w)}:(q,w)\rightarrow (q',w\alpha)\mid\alpha:q\rightarrow q'\in Q_1\}.\]
For an $\alpha\in \mathfrak{Q}_1$ define 
\[o(\alpha)=\left\{\begin{array}{l} 1\text{ if } \alpha\in Q_1\\-1 \text{ if }\alpha^{-1}\in Q_1\,.\end{array}\right.\]
The universal abelian covering quiver $\hat{Q}$, see Section $\ref{univcov}$, is given by the vertex set
\[\hat{Q}_0=Q_0\times\mathbb{Z}Q_1=\{(q,z_1)\mid q\in Q_0,z_1\in\mathbb{Z}Q_1\}\]
and the arrow set
\[\hat{Q}_1=\{(\alpha,z_1):(q,z_1)\rightarrow (q',z_1+e_{\alpha})\mid\alpha:q\rightarrow q'\in Q_1,\,z_1\in\mathbb{Z}Q_1\}.\]
The $k$-th universal abelian covering quiver is recursively defined by the vertex set
\[\hat{Q}^k_0=\hat{Q}_0^{k-1}\times\mathbb{Z}\hat{Q}_1^{k-1}=\{(q,z_1,\ldots,z_k)\mid q\in Q_0,z_l\in\mathbb{Z}\hat{Q}_1^{l-1}\}\]
and the arrow set
\begin{eqnarray*}\hat{Q}_1^{k}&=&\{(\alpha,z_1,\ldots,z_k):(q,z_1,\ldots,z_k)\rightarrow (q',z_1+e_{\alpha},\ldots,z_k+e_{(\alpha,z_1,\ldots,z_{k-1})})\\
&&\mid\alpha:q\rightarrow q'\in Q_1, z_l\in\mathbb{Z}\hat{Q}_1^{l-1}\text{ for }l=1,\ldots,k\}.\end{eqnarray*}
where we define $\hat{Q}^0=Q$. Note that $\hat{Q}_1^{k}=\hat{Q}_1^{k-1}\times\mathbb{Z}\hat{Q}_1^{k-1}$.

Fixing a vertex $q\in Q_0$ we can always consider those connected components  of $\tilde{Q}$ and $\hat{Q}^k$ such that the vertices $(q,1)$ and $(q,0,\ldots,0)$ are contained in these components. We again denote these subquivers by $\tilde{Q}$ and $\hat{Q}^k$. This is no restriction because we are only interested in stable and indecomposable representations respectively. Thus the support of such a representation is connected. Fix some vertex $(q',w)\in\tilde{Q}_0$ in this connected component. This means that $w=(q\mid\alpha_1\cdots\alpha_n\mid q')$ is a {\it path} in $Q_1$ and we may assume that $\alpha_i\neq\alpha_{i+1}^{-1}$ for all $i=1,\ldots,n-1$. We call such a path reduced in what follows. By $l(w)=n$ we denote the length of the path $w$ and, moreover, we define
\[h(i):=\frac{2i-1-o(\alpha_i)}{2}\]
and for $0\leq l\leq l(w)$ we define
\[c_1^l(w):=\sum_{i=1}^lo(\alpha_i)e_{\alpha_i}\in\mathbb{Z}Q_1\]
where $c_1^0(w)=0$. Furthermore, we recursively define
\[c_k^l(w):=\sum^l_{i=1}o(\alpha_i)e_{(\alpha_i,c^{h(i)}_1(w),\ldots,\ldots,c^{h(i)}_{k-1}(w))}\in\mathbb{Z}\hat{Q}_1^{k-1}\]
where again $c_k^0(w)=0$. Roughly speaking, $c_k^l(w)$ is the $k$-th coordinate of the vertex that we reach after $l$ steps in some universal abelian covering quiver, when walking along the path $w$. Now we can define a map $f_k:\tilde{Q}\rightarrow\hat{Q}^k$ by
\[f_k((q,w))=(q,c^{l(w)}_1(w),\ldots,c^{l(w)}_k(w))\] and
for some $\alpha:(q,w)\rightarrow (q',w\alpha)$ we define
\[f_k(\alpha)=(\alpha,(c^{l(w)}_i(w))_{i=1,\ldots,k}):(q,(c^{l(w)}_i(w))_{i=1,\ldots,k})\rightarrow (q',(c^{l(w\alpha)}_i(w\alpha))_{i=1,\ldots,k}).\]
Roughly speaking a path $w$ in $\tilde{Q}$ starting in $(q,1)$ is mapped to the same path in $\hat{Q}^k$. We get the image of such a path just by walking along ''the same arrows'' in $\hat{Q}^k$. But since the second quiver has cycles, different vertices and arrows can be mapped to the same vertices and arrows. But different paths are sent to different paths. Note that every arrow in both quivers corresponds to an arrow of the original quiver $Q$.

Obviously we have
\[c_k^{l(w\alpha)}(w\alpha)=c_k^{l(w)}(w)+e_{(\alpha,c^{l(w)}_1(w),\ldots,c^{l(w)}_{k-1}(w))}\]
for an arrow $\alpha\in Q_1$.
\begin{pro}
\begin{enumerate}
\item The maps $f_k$ are surjective for all $k$.
\item For $k\rightarrow\infty$ the map $f_k$ is injective.
\end{enumerate}
\end{pro}
{\it Proof.}
We first show that $f_k$ is surjective. As already mentioned we may consider the connected components such that $(q,1)\in\tilde{Q}_0$ and $(q,0,\ldots,0)\in\hat{Q}^k$. Thus let $(q',z_1,\ldots,z_k)\in \hat{Q}^k_0$ such that there exists a reduced path 
\[((q,0,\ldots,0)\mid\alpha_1\ldots\alpha_n\mid(q',z_1,\ldots,z_k))\]
in $\hat{Q}^{k}$ which corresponds to a reduced path $w=(q\mid\alpha_1\ldots\alpha_n\mid q')$ in $Q$.

We have $z_t=c_t^{l(w)}(w)$ and thus $f_k(q',w)=(q',z_1,\ldots,z_k)$.\\

Let $q\in Q_0$ and $w=\alpha_1\ldots\alpha_{n_k}\neq 1$ be a reduced path such that $w$ starts and terminates at $q$.  Moreover, assume that $q$ and $w$ are chosen such that $n_k$ is minimal satisfying \[f_k((q,1))=f_k((q,w))=(q,0,\ldots,0).\]
Then we claim that $f_{k+1}((q,w))\neq f_{k+1}(q,1)=(q,0,\ldots,0)$. Assume that this is not the case. Then we have
\[f_{k+1}((q,w))=(q,c^{l(w)}_1(w),\ldots,c^{l(w)}_{k+1}(w))=(q,0,\ldots,0)\] and, therefore, \[c_{t}^{l(w)}(w)=\sum^{l(w)}_{i=1}o(\alpha_i)e_{(\alpha_i,c^{h(i)}_1(w),\ldots,c^{h(i)}_{t-1}(w))}=0\] for all $t=1,\ldots,k+1$. Thus there exist $i,i'\in\{1,\ldots,n_k\}$ with $i\neq i'$ such that $c_l^{h(i)}(w)=c_l^{h(i')}(w)$ for all $l=1,\ldots,k$. But since $c_l^{h(i)}(w)-c_l^{h(i')}(w)=0$, this defines two vertices $(q',w_1)$ and $(q',w_1w_2)$ such that $f_k((q',w_1))=f_k((q',w_1w_2))$ and $l(w_2)<n_k$. But this contradicts the minimality of $n_k$. This already shows that $f_{\infty}$ is injective.\qed


Let $T_k:=(\mathbb{C}^{\ast})^{|\hat{Q}^{k-1}_1|}$.
Define \[M_d^{s}(Q)^{T,n}=(\ldots(M_d^s(Q)^{T_1})\ldots)^{T_n}.\]
Using Theorem $\ref{isom}$ we get the following:
\begin{satz}\label{bij}
For all dimension vectors $d$ there exists an $n\in\mathbb{N}_0$ such that we have \[M_d^{s}(Q)^{T,n'}\cong\bigcup_{\tilde{d}} M^s_{\tilde{d}}(\tilde{Q})\] for all $n'\geq n$ where $\tilde{d}$ ranges over all
equivalence classes that are compatible with $d$.
\end{satz}

Concerning the Euler characteristic of quiver moduli spaces we get the following corollary:
\begin{kor} Let $Q$ be a quiver with dimension vector $d$. Then for the Euler characteristic of the moduli space $M^s_d(Q)$ we have
\[\chi(M^s_d(Q))=\sum_{\tilde{d}}\chi(M^s_{\tilde{d}}(\tilde{Q})),\]
where $\tilde{d}$
ranges over all equivalence classes being compatible with $d$.
\end{kor}

Thus, if we are interested in the Euler characteristic, we may always assume that torus fixed points are given as representations of the universal covering quiver which has no cycles. For this quiver the representations theory often simplifies in comparison to the universal abelian covering quiver. 
\begin{bem}
\end{bem}
\begin{enumerate}
\item Every connected component of the universal abelian covering quiver of the Kronecker quiver $K(m)$ is an infinite bipartite $(m-1)$-dimensional honeycomb lattice. This means that we have given an orientation such that every vertex is either a source or a sink.
\item Every connected component of the universal covering quiver of the Kronecker quiver $K(m)$ is a infinite bipartite regular $m$-tree. 
\end{enumerate}

We end up this section with the following definition:
\begin{defi}
Let $Q$ be a quiver and $\Theta\in\mathrm{Hom}(\mathbb{Z}Q_0,\mathbb{Z})$ a linear form. A tuple consisting of a finite subquiver $\mathcal{Q}$ of $\hat{Q}$ (resp. $\tilde{Q})$ and a dimension vector $d\in\mathbb{N}\mathcal{Q}_0$ such that $M^{s}_{d}(\mathcal{Q})\neq\emptyset$, where we consider the stability induced by $\Theta$, is called localization data.
\end{defi}
There exists an equivalence relation on the set of localization data obtained by translating the vertices by $\chi\in\mathbb{Z}Q_1$. In the following, we will always consider localization data up to this equivalence.

A localization data always comes along with an embedding into some covering quiver. This induces a colouring of the arrows $c:\mathcal{Q}_1\rightarrow Q_1$. If we forget about this colouring we call such a data uncoloured localization data. The purpose of it is that fixing an uncoloured localization data there can exist many colourings that induce different localization data.

\section{Localization in Kronecker moduli spaces}\label{lockro}
\noindent Since the main focus of the paper is on the generalized Kronecker quiver, in this section we apply the introduced machinery to this case. We first recall some properties of Kronecker moduli spaces. Then we investigate stable bipartite quivers which representations are torus fixed points of these moduli spaces when colouring the arrows appropriately. In particular, we construct a class of localization data which grows exponentially with the dimension vector.
\subsection{Kronecker moduli spaces}
\noindent Let $K(m)$ be the generalized Kronecker quiver having two vertices and $m$ arrows between them, i.e.:
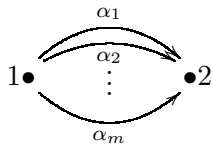
\[
\begin{xy}
\xymatrix@R10pt@C20pt{
\\1\bullet\ar@/_1.5pc/[rr]_{\alpha_m}\ar@/^1.5pc/[rr]^{\alpha_1}\ar@/^1.0pc/[rr]_{\alpha_2}&\vdots &\bullet 2\\&&
}
\end{xy}\]
A representation of this quiver with dimension vector $(d,e)$ is given by two $\mathbb{C}$-vector spaces $V$ and $W$ of dimensions $d$ and $e$ and an $m$-tuple of linear maps 
\[(X_1,\ldots,X_m)\in\bigoplus_{i=1}^m \mathrm{Hom}(V,W)=R_{d,e}(K(m)).\]
The group $(Gl(V)\times Gl(W))$ acts on $R_{d,e}(K(m))$ via simultaneous base change. Since the scalar matrices act trivially, the group action factorises through the quotient $(Gl(V)\times
Gl(W))/\mathbb{C}^{\ast}$. For $\Theta=(1,0)$ the slope function
$\mu:\mathbb{N}^2\backslash\{0\}\rightarrow\mathbb{Q}$ is given by
\[\mu(d,e):=\frac{d}{d+e}.\]
Thus we obtain the following criterion for the (semi-)stability of Kronecker representations:
\begin{lem}
A point $(X_1,\ldots,X_m)\in
R_{d,e}(K(m))$ is semistable (resp. stable) if and only if for all proper subspaces $0\neq U\subsetneq V$ the following holds:
\[\dim\sum\limits_{k=1}^m X_k(U)\geq  \frac{e}{d}\cdot\dim U
\text{ }(\text{resp. }\dim\sum\limits_{k=1}^m X_k(U)> \frac{e}{d}\cdot\dim U ).\]
\end{lem}

In particular, a dimension vector is $\Theta$-indivisible if and only if $d$ and $e$ are coprime. In the following, we call the geometric quotient $M_{d,e}^m:=M^s_{d,e}(K(m))$ Kronecker moduli space. Using standard methods from Algebraic Geometry, see for instance \cite{sha} and \cite{sha2}, we get by use of Theorem \ref{kingsatz} the following:
\begin{kor}
Let $d,e\in\mathbb{N}$ such that $\gcd(d,e)=1$.
The corresponding Kronecker moduli space $M^m_{d,e}$ is a compact complex manifold. Furthermore, there exists a continuous map
\begin{center}
$\Pi:R_{d,e}^s(K(m))\rightarrow M^m_{d,e}$ 
\end{center}
such that the $\Pi$-fibres are exactly the orbits under the group action.
\end{kor}
\begin{bem}
\end{bem}
\begin{enumerate}
\item Note that for $m=1$ there exist only indecomposable (resp. stable) representations of dimensions $(1,0),\,(0,1)$ and $(1,1)$ and for $m=2$ the only cases of interest are the dimension vectors $(d,d)$, $(d,d+1)$ and $(d+1,d)$ for $d\in\mathbb{N}$. The roots $(d,d)$ for $d\geq 2$ are no Schur roots which means that the moduli spaces $M_{d,d}^2$ are empty.

Furthermore, $(d,d+1)$ is a real Schur root which means that the moduli space is a point. Thus we will assume that $m\geq 3$ if we do not explicitly say anything else.
\end{enumerate}
We state some helpful properties of Kronecker moduli spaces:
\begin{pro}\label{propofkron}
\begin{enumerate}
\item There exist isomorphisms of moduli spaces $M^m_{d,e}\simeq M^m_{e,d}$ and
$M^m_{d,e}\simeq M^m_{me-d,e}$. \item The dimension of the moduli spaces is given by \[\dim\,M^m_{d,e}=1-d^2-e^2+dem\] if
$M^m_{d,e}\neq \emptyset.$ \item We have $M^m_{d,e}\neq\{pt\}$ if and only if \[\frac{\textstyle m-\sqrt{m^2-4}}{\textstyle
2}<\frac{\textstyle e}{\textstyle d}<\frac{\textstyle
m+\sqrt{m^2-4}}{\textstyle 2}\text{ holds.}\]
\end{enumerate}
\end{pro}
{\it Proof.}
We obtain the first isomorphism by considering the map
\[\overline{(X_1,...,X_m)}\rightarrow\overline{(X_1^T,...,X_m^T)}.\]
The second one is obtained via the reflection functor, see Theorem $\ref{kspi}$.
The second part is a special case of the fifth part of Remark $\ref{bem1}$.

If
$M^m_{d,e}\neq\{pt\}$ holds, \[\frac{\textstyle
m-\sqrt{m^2-4}}{\textstyle 2}\leq\frac{\textstyle e}{\textstyle
d}\leq\frac{\textstyle m+\sqrt{m^2-4}}{\textstyle 2}\] follows from the second part of the proposition. But $K(2)$ with dimension vector $(d,d)$ is the only case such that equality holds. 

If the inequalities are satisfied properly, $(d,e)$ is an imaginary Schur root, see \cite{kac}. In particular, we have $\langle (d,e),(d,e)\rangle=d^2+e^2-dem\leq 0$. Following \cite{sch} there exists an open subset of $R_{d,e}(K(m))$ which contains those representations which are stable in the sense of King with $\tilde{\Theta}((d',e')):=\langle(d',e'),(d,e)\rangle-\langle (d,e),(d',e')\rangle$.\qed

\subsection{Localization data of the Kronecker quiver}
\noindent In this subsection we investigate the support of the dimension vectors which arise from the localization method in detail. Moreover, we investigate stable bipartite quivers and the possibilities of colouring their arrows so that stable representations of such quivers become torus fixed points.\\

Let $(d,e)\in\mathbb{N}^2$ be a dimension vector of the Kronecker quiver and let
\[X=((V,W),(X_1,\ldots,X_m))\in (M^m_{d,e})^T\] be a fixed point. From the considerations of the third section we get a morphism of algebraic groups $\varphi:T\rightarrow (Gl(V)\times Gl(W))/\mathbb{C}^{\ast}$, for which we can choose a lift
$\psi:T\mapsto Gl(V)\times Gl(W)$. It can be decomposed into two morphisms $\psi_1:T\mapsto Gl(V)$ and $\psi_2:T\mapsto
Gl(W)$.

Let 
\[V=\bigoplus_{\chi\in X(T)}V_{\chi} \text{ and } W=\bigoplus_{\chi\in
X(T)}W_{\chi}\] 
be the weight space decompositions with respect to $\psi_1$ and
$\psi_2$ respectively. They satisfy
\[X_k(V_{\chi})\subseteq W_{\chi+e_k}\]
for all $\chi\in X(T)\cong\mathbb{Z}^m$ and $k=1,\ldots ,m$.\\

The universal abelian covering quiver
$\hat{K}(m)$ has vertices $(1,\chi)$ and $(2,\chi)$, where $\chi$ runs through all characters of $X(T)$, and arrows \[(1,\chi)\rightarrow (2,\chi+e_k)\] for every $k\in\{1,\ldots,m\}$ and every $\chi\in\mathbb{Z}^m$.

For every fixed point there exists a unique dimension vector $\hat{d}$ given by
\[\hat{d}_{1,\chi}=\dim V_{\chi}\text{ and
}\hat{d}_{2,\chi}=\dim W_{\chi}\] for $(1,\chi),(2,\chi)\in\hat{K}(m)_0.$

The other way around consider $\hat{K}(m)$ and a dimension vector $\hat{d}\in\mathbb{N}\hat{K}(m)_0$. A stable representation of this quiver corresponds to a torus fixed point with dimension vector $(d,e)$ where
\[d=\sum_{\chi\in X(T)}\hat{d}_{1,\chi} \text{ and }e=\sum_{\chi\in X(T)}\hat{d}_{2,\chi}.\]
In what follows, we call $(d,e)$ dimension type of the representation.
\begin{defi}
Let $Q$ be a quiver with a fixed linear form $\Theta\in\mathrm{Hom}(\mathbb{Z}Q_0,\mathbb{Z})$. A tuple consisting of the quiver $Q$ and a dimension vector $d\in\mathbb{N}\mathcal{Q}_0$ is called stable if $M^s_d(Q)\neq\emptyset$ where we consider the stability induced by $\Theta$.
\end{defi}
If it is clear which dimension vector we consider, we will simply call such a tuple stable quiver.
\begin{bem}
\end{bem}
\begin{enumerate}
\item
The stability condition for representations of $\hat{K}(m)$ is induced by the original linear form $\Theta=(1,0)$. It is given by
\[\mu(\hat{d})=\frac{\sum_{\chi\in X(T)}\hat{d}_{1,\chi}}{\sum_{\chi\in X(T)}\hat{d}_{1,\chi}+\hat{d}_{2,\chi}}.\]
\end{enumerate}
Let $\hat{K}(m)_c$ be a connected component of $\hat{K}(m)$. It is bipartite  and, moreover, there exists an embedding
$\lambda:(\hat{K}(m)_{c})_0\rightarrow\mathbb{Z}^m$ such that $\lambda(q,\chi)=\chi$ for $q=1,2$. In the following let $I\cup J$ be the decomposition of the vertex set into sources and sinks. We may assume that they are elements of $\mathbb{Z}^m$. Let $R\subset I\times J$ be the set of arrows. Then we have
\[(i,j)\in R\Leftrightarrow j=i+e_k\] for
$k\in\{1,\ldots,m\}$. This defines a map $c:R\rightarrow\{1,\ldots,m\}$, which we call colouring in the following, by setting $c(i,j)=k$ if $j=i+e_k$. Obviously, the set $R$ and the colouring $c$ are already uniquely determined by the vertex set $I\cup J$. Nevertheless, they play an important role because they will describe different localization data for a fixed uncoloured localization data or a fixed stable bipartite quiver.

In conclusion, a fixed point $X$ determines a tuple $(I,J,\hat{d})$ which is unique up to translation by a vector $\mu\in\mathbb{Z}^m$. In what follows we always consider such tuples up to translation.

For a bipartite quiver with vertex set $I\cup J$ and fixed dimension vector $d$ define the sets
\[A_i:=\{j\in J\mid \alpha:i\rightarrow j\in Q_1,d_j\geq 1\}\text{ and }A_j:=\{i\in I\mid \alpha:i\rightarrow j\in Q_1,d_i\geq 1\}.\]
Furthermore, define $A_{I}=\bigcup_{i\in I}A_i$ and $R_i=|A_i|$ and $R_j=|A_j|$. 
\begin{defi}
A bipartite quiver is called $m$-bipartite if we have for all sources $i\in I$ and all sinks $j\in J$ that $R_i,R_j\leq m$.
\end{defi}
\begin{bem}\label{anzko}
\end{bem}
\begin{enumerate} \item Consider a stable bipartite quiver $(\mathcal{Q},\hat{d})$ without oriented and unoriented cycles of dimension type $(d,e)$ with $\mathcal{Q}=(I\cup J,R)$ such that there exists at most one arrow between every two vertices. Choose a colouring of the arrows $c:\mathcal{Q}_1\rightarrow\{1,\ldots,m\}$. Then we get a localization  data if $\mathcal{Q}$ and $c$ satisfy the following conditions:
\begin{enumerate}
\item The quiver $\mathcal{Q}$ is $m$-bipartite.
\item For all $(i,j),(i,j')\in R$ such that $j\neq j'$ we have
$c(i,j')\neq c(i,j)$. \item Analogously, for $(i,j),(i',j)\in
R$ such that $i\neq i'$ we have $c(i,j)\neq c(i',j)$.
\end{enumerate}
We call a colouring satisfying these conditions stable.
\item It is also easy to check that if $i\in I$ is a source such that $\hat{d}_i=1$, then we have
$m\geq R_i>\frac{e}{d}$.\end{enumerate}

Note that, fixing a vertex $q\in\mathcal{Q}_0$ and setting $c(q):=0$ every colouring of the arrows induces a colouring of the vertices $c:\mathcal{Q}_0\rightarrow \mathbb{Z}^m$. 

\begin{bem}\label{darstbip}
\end{bem}
\begin{enumerate}
\item Because of Lemma \ref{zusa} the quiver of a localization data has to be connected. Otherwise there would exist an exact sequence contradicting the stability condition.
\item In order to test an $m$-bipartite quiver for stability, we do not need to consider an explicit representation. We can rather consider an arbitrary representation $X$ of this dimension satisfying for all $j\in J$ and all subsets $A'_j\subseteq A_j$ with $R'_j:=|A'_j|$ the following property:
\[\dim\bigcap_{i\in A'_j}X_{\alpha}(X_i)=\max\{0,\sum_{i\in A'_j}\dim X_{\alpha}(X_i)-(R'_j-1)\dim X_j\}.\]
Indeed, consider a bipartite quiver of the form
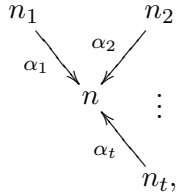
\[
\begin{xy}
\xymatrix@R15pt@C10pt{ n_1\ar[rd]_{\alpha_1}&&n_2\ar[ld]_{\alpha_2}\\
&n&\vdots\\
&&n_t\ar[lu]^{\alpha_t},}
\end{xy}
\]
where $(n,n_1,\ldots,n_t)$ denotes the dimension vector. If $n_i\leq n$ for all $1\leq i\leq t$, there always exists a representation $X$ of this quiver such that for all tuples of linear maps $X_{\alpha_{i_1}},\ldots,X_{\alpha_{i_k}}$ with $1\leq k\leq t$ and $1\leq i_1<i_2<\ldots<i_k\leq t$ the dimension of the intersections of the images is minimal. One verifies the existence and the dimension formula by induction on the number of arrows.
\end{enumerate}
\begin{bem}\label{indPfeile}
\end{bem}
\begin{enumerate}
\item Fixing a stable $m$-bipartite quiver without cycles it may happen that different colourings of the arrows lead to different types of localization data. For instance, if we consider a colouring $c$ such that this colouring induces a weight space of weight $\chi$ and one of weight  $\chi-e_k$, we have an arrow $\alpha:\chi-e_k\rightarrow \chi$ and also a linear map
\[X_{\alpha:\chi-e_k\rightarrow\chi}:V_{\chi-e_k}\rightarrow V_{\chi}.\]
We call such an arrow induced. Obviously, the dimension of the corresponding moduli space of the universal abelian cover increases at least by one in comparison to the dimension of the moduli space corresponding to the bipartite quiver.

\hspace{0.3cm}Moreover, it can happen that two different vertices are identified when choosing a colouring. Thus a colouring can induce two types of cycles. But Theorem $\ref{bij}$ put things right. In particular, after suitable many localization steps the remaining torus fixed points are representations of the universal covering quivers which has no cycles. 
\end{enumerate}
\begin{lem}\label{zykel}
Let $(\mathcal{Q},d)$ be a stable $m$-bipartite quiver without cycles and $c,c'$ stable colourings of the arrows. Then we have:
\begin{enumerate}
\item By colouring the arrows with $c$ we obtain a localization data. 
\item Fix $c$ and $c'$ such that $c$ induces no cycles and $c'$ induces at least one cycle. Moreover, let $\dim(M_{\mathcal{Q}},c)$ and $\dim(M_{\mathcal{Q}},c')$ be the dimensions of the resulting moduli spaces. We have
\[\dim(M_{\mathcal{Q}},c)\leq\dim(M_{\mathcal{Q}},c').\]
\end{enumerate}
\end{lem}
{\it Proof.}
Fixing a stable $m$-bipartite quiver and a stable colouring of the arrows we obtain a localization data. Every stable representation of $\mathcal{Q}$ induces a stable representation of $\hat{K}(m)$, no matter if the colouring leads to cycles or not. Induced arrows let the dimension of the moduli space increase. Thus it remains to prove that the dimension of the moduli space increases if a colouring induces a cycle which does not come from an induced cycle.
Let $(\mathcal{Q},c)$ and $(\mathcal{Q},c')$ be the two resulting subquivers of $\hat{K}(m)$ and $d(c)$ and $d(c')$ be the resulting dimension vectors respectively. Assume $c'(j_1)=c'(j_2)$ and let $j_{1,2}$ be the corresponding vertex of $(\mathcal{Q},c')$, i.e. $d(c')_{j_{1,2}}=d_{j_1}+d_{j_2}$. We have $R_{j_1},R_{j_2}\geq 1$ in $Q$. Define $\dim A_j=\sum_{i\in A_j}d_i$. 
Then we have for the colouring $c'$ producing this cycle
\begin{eqnarray*}
\dim(M_{\mathcal{Q}},c')&=&\dim(M_{\mathcal{Q}},c)+d_{j_1}^2+d_{j_2}^2-(d_{j_1}+d_{j_2})^2-\dim A_{j_1}d_{j_1}\\&&-\dim A_{j_2}d_{j_2}+(\dim A_{j_1}+\dim A_{j_2})(d_{j_1}+d_{j_2})\\
&=&\dim(M_{\mathcal{Q}},c)-2d_{j_1}d_{j_2}+\dim A_{j_1}d_{j_2}+\dim A_{j_2}d_{j_1}\geq \dim(M_{\mathcal{Q}},c).
\end{eqnarray*}
Indeed, we have $\dim A_{j_k}\geq d_{j_k}$ for $k=1,2$ because of the stability of $Q$. The case $c'(i_1)=c'(i_2)$ is proved in the same way.\qed
\begin{defi}\label{verk}
A localization data $(\mathcal{Q},\hat{d})$ is called localization data of type one if $\hat{d}_q\in\{0,1\}$ for all $q\in \mathcal{Q}_0$.
\end{defi}
  
\subsection{Stability of bipartite quivers}\label{stabilitaet}
\noindent In this section we study how to construct new stable bipartite quivers by glueing stable bipartite quivers of smaller dimension types. Hence each stable colouring gives rise to some localization data. This gives a huge class of such quivers for every fixed dimension type. Note that this glueing method is also used to construct stable tree modules of the Kronecker quiver, see \cite{wei3}.\\

Let $Q=(I\cup J,Q_1)$ and $Q'=(I'\cup J',Q_1')$ be two bipartite quivers with $j\in J$, $j'\in J'$. Define the bipartite quiver \[Q_{j,j'}(Q,Q')=(I\cup I'\cup J\backslash j \cup J'\backslash j'\cup j'' ,Q_1'')\] 
such that $\alpha:i\rightarrow j_1\in Q_1''$ if and only if $\alpha:i\rightarrow j_1\in Q_1$ or $\alpha:i\rightarrow j_1\in Q_1'$ with $j_1\neq j,j'$ and  $\alpha:i\rightarrow j''\in Q_1''$ if and only if  $\alpha:i\rightarrow j\in Q_1$ or $\alpha:i\rightarrow j'\in Q_1'$.

Thus the new quiver is generated by the former ones by identifying two vertices of the set of sinks of these quivers.
\begin{defi}
The quiver $Q_{j,j'}(Q,Q')$ is called the glueing quiver of $Q$ and $Q'$ and the vertices $j,j'=j''$ the glueing vertices.
\end{defi}
\begin{defi}
Let $(Q,d)$ be a tuple consisting of a bipartite quiver with sources $I$ and $d\in\mathbb{N}Q_0$ a dimension vector. A subquiver of $Q$ with sources $I'$ is called boundary quiver if there exists one $i_0\in I'$ such that $|A_{i_0}\cap A_{I\backslash I'}|=1$ and $|A_i\cap A_{I\backslash I'}|=0$ for all $i\in I'$ with $i\neq i_0$. A boundary quiver is called proper boundary quiver if it does not contain any other boundary quiver.
\end{defi}

Note that if $d_q\geq 1$ for all $q\in Q_0$ this means that boundary quivers are such subquivers which only have one common sink with the remainder of the quiver.

Fixing a representation $X$ of $Q$, we abbreviate the dimension of the image of a subspace $U=\oplus_{i\in I} U_i$ with $U_i\subset X_i$ to $d_U$. Explicitly, we define
\[d_U:=\dim\sum_{\substack{i\in I\\\alpha\in\{\beta\mid t(\beta)=i\}}}X_{\alpha}(U_i).\]

Fixing a coprime dimension vector $(d,e)$ with $d\geq 1$ we now determine a unique dimension vector $(d_s,e_s)$ such that we are able to construct new stable bipartite quivers of dimension type $(d_s+(k+l)d,e_s+(k+l)e)$ by glueing quivers of the types $(d_s+kd,e_s+ke)$ and $l(d,e)$.

Fixing some dimension vector $(d,e)$, we first show that there exists a coprime dimension vector $(d_s,e_s)$ such that $d_s\leq d$ and $e_s\leq e$ satisfying the conditions
\begin{enumerate}
\item $\frac{e+e_s}{d+d_s}d<e+1$ if $d\neq 1$.
\item $\frac{e+e_s}{d+d_s}d>e$ if $d\neq 1$.
\item $\frac{e_s-1}{d_s}< \frac{e}{d}$ if $d\neq 1$ and $(e_s-1)d=ed_s$ if $d=1$.
\item $\frac{e+e_s}{d+d_s}d'<\lceil\frac{e}{d}d'\rceil\hspace{0.3cm}\forall\, 1\leq d'<d$.
\item $\gcd(d+d_s,e+e_s)=1.$
\end{enumerate}
These are conditions which should intuitively be satisfied in order to be able to glue stable quivers of dimension types $(d_s,e_s)$ and $(d,e)$ to get one of dimension type $(d_s+d,e_s+e)$. We refer to these conditions as glueing conditions. We will see that these conditions are also sufficient. The first property is equivalent to the following:
\begin{eqnarray*}
de+de_s< de+d+d_se+d_s &\Leftrightarrow&  de_s < d+d_se+d_s\\
\Leftrightarrow  d(e_s-1)< d_s(e+1)&\Leftrightarrow&  \frac{e_s-1}{d_s} <\frac{e+1}{d}.
\end{eqnarray*} 
The second one is equivalent to:
\[
ed+e_sd > ed+ed_s\Leftrightarrow \frac{e_s}{d_s}>\frac{e}{d}.
\]
Therefore, it suffices to verify the second and third property because the first one follows from the third one.
\begin{lem}\label{startexist}
Let $(d,e)\in\mathbb{N}^2$ such that $d\leq e$ and $d,e$ are coprime. There exists a coprime dimension vector $(d_s,e_s)$ satisfying the glueing conditions. It is uniquely determined if we also assume that $d_s\leq d$ and $e_s\leq e$. 
\end{lem}
{\it Proof.}
We first consider the special case $d=1$. It is easy to see that $(0,1)$ satisfies these properties for $(d,e)=(1,n)$ with $n\in\mathbb{N}$.\\

If $d\geq 2$, we already have $e\geq 3$. Choose $d_s\in\mathbb{N}$ minimal such that $d\mid 1+ed_s$. This is possible because $\gcd(d,e)=1$ and, therefore, there exist $\lambda,\mu\in\mathbb{Z}$ such that 
\[\lambda d=1-\mu e.\]
Define
\[e_s=\frac{1+d_se}{d}.\]
Because of the choice of $d_s$, we have $e_s\in\mathbb{N}$.

Moreover, we get
\[-e(d+d_s)+d(e+e_s)=-ed-ed_s+de+d_se+1=1.\]
It follows that $\gcd(d+d_s,e+e_s)=1$.

Now we get
\[\frac{e_s}{d_s}=\frac{1+d_se}{dd_s}>\frac{e}{d}\]
and also
\[\frac{e_s-1}{d_s}=\frac{d_se-d+1}{dd_s}<\frac{e}{d}.\]
Thus it remains to prove the fourth property. By an easy calculation we get
\[\frac{e+e_s}{d+d_s}=\frac{e}{d}\left(\frac{ed+ed_s+1}{ed+ed_s}\right)=\frac{e}{d}\left( 1+\frac{1}{ed+ed_s}\right) .\]
Moreover, since 
\[\lceil\frac{e}{d}d'\rceil-\frac{e}{d}d'\geq\frac{1}{d}\]
 and 
\[\frac{d'}{d(d+d_s)}<\frac{1}{d+d_s}\]
for each $d'<d$, the existence of such a vector follows.\\

If $(d'_s,e'_s)$ is another dimension vector satisfying the desired properties, it is straightforward that the glueing conditions imply
\[d'_se_s-\frac{d'_s}{d}-1<d_se'_s-1<d'_se_s.\]
But since $\frac{d'_s}{d}<1$ it follows that $d'_se_s=d_se'_s$. But since $(d_s,e_s)$ and $(d'_s,e'_s)$ are both coprime we already have $(d_s,e_s)=(d'_s,e'_s)$.\qed

In what follows, we call a vector $(d_s,e_s)$ satisfying these properties starting vector for $(d,e)$. In the remainder of the section we assume that $(d,e)$ is coprime and $(d_s,e_s)$ is the corresponding starting vector as constructed in Lemma $\ref{startexist}$.
\newpage
\begin{bem}\label{darstellung}
\end{bem}
\begin{enumerate}
\item If we want to decompose a coprime dimension vector $(d,e)$ into
\[(d,e)=(d_s,e_s)+k(d',e')\]
such that $(d',e')$ and $(d_s,e_s)$ satisfy the glueing conditions, we can proceed as follows: let $e'\in\mathbb{N}$ minimal such that
\[e\mid 1+de'\text{ and }d'=\frac{1+e'd}{e}.\] Now we compute $d_s$ and $e_s$ as before. It can be seen easily that these numbers satisfy the glueing conditions. Indeed, one checks that
\[\frac{e-e_s}{e'}=\frac{d-d_s}{d'}.\]
It follows that $e'\mid e-e_s$ and $d'\mid d-d_s$ 
because $\gcd(d',e')=1$ and, trivially, $e-e_s,d-d_s\in\mathbb{N}$ hold. Now define $k=\frac{d-d_s}{d'}$. 
\end{enumerate}
We need other properties of these natural numbers. By use of $e_sd-ed_s=1$ we get
\[(ke+e_s)(k'd+d_s)+k-k'=(kd+d_s)(k'e+e_s)\]
where $k,\,k'\in\mathbb{N}$. For $d_1=k'd+d'\in\mathbb{N}$ with $0\leq d'< d$ and $0<d_1\leq kd+d_s$ define a map
\[f(d_1)=\min\{n\in\mathbb{N}\mid \frac{(ke+e_s)d_1+n}{kd+d_s}\in\mathbb{N}\}.\]
Note that $f$ is injective because $\gcd(kd+d_s,ke+e_s)=1$.
Then we get the following lemma:
\begin{lem}
Let $d_s$, $e_s$, $d$, $e$ fulfil the glueing conditions. Then we have
\[(ke+e_s)(k'd+d_s)+k-k'=0\modu (kd+d_s)\]
for all $k'\leq k$.

Let $d_1=k'd+d'$ with $0\leq d'<d$. In particular, we have $f(d_1)=k-k'$ if $d'=d_s$ and thus $f(d_1)\geq k+1$ if $d'\neq d_s$. 
\end{lem}

Now we show how to get a stable bipartite quiver of dimension type $(d_s+(k+l)d,e_s+(k+l)e)$ by glueing a stable bipartite quiver of type $(d_s+kd,e_s+ke)$ and certain quivers of type $(ld,le+1)$. We again point out Remark $\ref{darstbip}$. Thus we do not consider specific representations, but those satisfying the properties mentioned in the remark.\\

In the following, if we fix a bipartite quiver $\mathcal{Q}$, we always additionally fix a dimension vector $\hat{d}\in\mathbb{N}\mathcal{Q}_0$. In abuse of notation we do not always mention it and, moreover, if we glue two bipartite quivers the dimension vector of the glueing quiver is denoted by $\hat{d}$ again. We just additionally specify the dimension corresponding to the glueing vertex. The remaining vertices keep the dimension.
 
Let $\mathcal{S}_{ld,le+1}^m$ be the set of tuples consisting of a connected $m$-bipartite quiver $\mathcal{Q}$ of dimension type $(ld,le+1)$ and a sink $j$ with $\hat{d}_j\geq 1$ satisfying the following properties:
\begin{enumerate}
\item There exists a representation $T$ (with the corresponding dimension vector) of the quiver such that for every $d'$-dimensional subspace $U$ we have
\[d_U>\frac{(k+l)e+e_s}{(k+l)d+d_s}d'.\]
\item After decreasing the dimension of the sink $j$ by one, the resulting quiver is connected and the corresponding factor representation of $T$ is semistable, i.e. the quiver is semistable.
\end{enumerate}
Let $\mathcal{T}_{d,e}^m$ be the set of all stable $m$-bipartite quivers of dimension type $(d,e)$.
\begin{satz}\label{verkl}
Let $d$ and $e$ be coprime, $d,d_s,e,e_s$ fulfil the glueing conditions and let $k\in\mathbb{N}$. Let $T^0\in\mathcal{T}^m_{d_s+kd,e_s+ke}$ and $(T^1,j_1)\in\mathcal{S}_{ld,le+1}^m$. Moreover, let $j_0$ with $\hat{d}_{j_0}$ be a sink of $T_0^0$ such that $R_{j_0}+R_{j_1}\leq m$. Then $Q_{j_0,j_1}(T^0,T^1)$ with glueing vertex $j_2$ where $\hat{d}_{j_2}:=\hat{d}_{j_0}+\hat{d}_{j_1}-1$ is an element of $\mathcal{T}^m_{d_s+(k+l)d,e_s+(k+l)e}$.
\end{satz}
{\it Proof.} 
For some subspace $U$ of one of the two subquivers we denote by $d_U$ the dimension of its image corresponding to its original quiver and by $d'_U$ the dimension of its image corresponding to the glueing quiver.

Given a stable representation $S$ of $T^0$ and a representation $T$ of $T^1$ satisfying the condition from above we consider the following representation of $Q_{j_0,j_1}(T^0,T^1)$: the corresponding semistable factor representation of $T$ induces a one-dimensional subspace of $T_{j_1}$ which we identify with an arbitrary one-dimensional subspace of $S_{j_0}$.

First let $U$ be a $d'$-dimensional subspace corresponding to $T$ such that $d'<ld$. Then by definition we have \[\frac{(k+l)e+e_s}{(k+l)d+d_s}d'<d_U=d'_U.\]
If $d'=ld$, the same inequality follows from $d_U=le+1$ which follows from the glueing condition.

Since we also have
\begin{equation}\label{eq1}\frac{e_s+ke}{d_s+kd}>\frac{e_s+(k+l)e}{d_s+(k+l)d},\end{equation} see the properties of the dimension vectors, the same follows for subspaces of $S$.

It remains to prove that subspaces composed of subspaces of both subquivers fulfil the stability condition. Thus let $U'$ and $U''$ be two subspaces of dimensions $1\leq d'\leq ld$ and $1\leq d''\leq kd+d_s$ respectively such that we have proper inequality at least once. Here $U'$ corresponds to $T$ and $U''$ to $S$.

Now it suffices to prove that
\[d'_{U'\oplus U''}\geq d_{U'}+d_{U''}-1\geq\frac{le}{ld}d'+d_{U''}>\frac{(k+l)e+e_s}{(k+l)d+d_s}(d'+d'')\]
where the first inequality follows by construction and the second inequality follows from the semistability of the quiver obtained from $T^1$ after decreasing the dimension of the vertex $j_1$ by one. The last inequality is equivalent to
\[d_{U''}>\frac{(k+l)e+e_s}{(k+l)d+d_s}d''+\frac{d'}{d((k+l)d+d_s)}\]
using $e_sd-d_se=1$.

By the preceding lemma together with the assumption we have
\[d_{U''}\geq \frac{(ke+e_s)d''+f(d'')}{kd+d_s}.\]
First let $d''<kd+ds$. Assuming without loss of generality that $d'=ld$, it remains to prove that
\[ld''+((k+l)d+d_s)f(d'')>l(kd+d_s).\]
But this is easily verified.

Finally, let $d''=kd+d_s$ and $d'=l'd+d_1<ld$ with $0\leq d_1<d$.
We have
\[\frac{(k+l)e+e_s}{(k+l)d+d_s}(kd+d_s)=ke+e_s-\frac{l}{(k+l)d+d_s}\]
again using $e_sd-ed_s=1$. Thus it remains to prove 
\[\lceil\frac{e}{d}(l'd+d_1)\rceil=l'e+\lceil\frac{ed_1}{d}\rceil>\frac{(k+l)e+e_s}{(k+l)d+d_s}(l'd+d_1)-\frac{l}{(k+l)d+d_s}\]
what follows from the fourth glueing condition and inequality (\ref{eq1}) together with $l>l'$.\qed

If $T^0$ and $T^1$ satisfy the condition of the theorem we call $T^0$ starting quiver for $T^1$. Now we apply the result to specific quivers. Therefore, let $T\in \mathcal{T}^m_{d,e}$. Starting with this quiver, we construct new quivers $\hat{T}$ of dimension type $(d,e+1)$ in one of the following ways:
\begin{enumerate} 
\item Choose an $i\in I$ such that $R_i<m$ and define the new quiver by the vertex set $\hat{T}_0=T_0\cup\{ j\}$ and the arrow set $\hat{T}_1=T_1\cup\{\alpha:i\rightarrow j\}$. Finally, let $\hat{d}_j=1$.\item Choose a vertex $j\in J$ with $1<R_j<m$ and increase the dimension of the vertex by one.
\item Choose a vertex $j\in J$ such that \[\hat{d}_j<\sum_{i\in A_{j}} \hat{d}_i\] 
and increase the dimension of the vertex $j$ by one.
\end{enumerate}
Denote the set of the resulting quivers by $\hat{\mathcal{T}}_{d,e}^m$ and refer to $j$ as modified vertex. Given a representation $X$ of $T\in\mathcal{T}_{d,e}^m$ we can modify it under consideration of Remark \ref{darstbip} in order to get a representation of $\hat{T}$.
\begin{kor}\label{korverkl}
Let $d,d_s,e,e_s$ be as before and let $k\in\mathbb{N}$. Moreover, let $T^0\in\mathcal{T}^m_{d_s+kd,e_s+ke}$ and $T^1\in\hat{\mathcal{T}}^m_{d,e}$ with modified vertex $j_1$. Further let $j_0$ with $\hat{d}_{j_0}\geq 1$ be a sink of $T_0^0$ such that $R_{j_0}+R_{j_1}\leq m$. Then $Q_{j_0,j_1}(T^0,T^1)$ with glueing vertex $j$, where $\hat{d}_j:=\hat{d}_{j_0}+\hat{d}_{j_1}-1$, is an element of $\mathcal{T}_{d_s+(k+1)d,e_s+(k+1)e}$.
\end{kor}
{\it Proof.}
We just have to check the two conditions stated before the preceding theorem. Thus let $U$ be a $d'$-dimensional subspace of a modified representation $\hat{X}$ of $T^1$. Since $T^1$ results from a stable quiver we have $d_U>\frac{e}{d}d'$.

Moreover, by the fourth glueing condition it follows that \[\frac{(k+1)e+e_s}{(k+1)d+d_s}d'<\lceil\frac{e}{d}d'\rceil\leq d_U.\]
If $d'=d$, the same inequality follows from the first property together with \[\hat{d}_{j_1}\leq\sum_{i\in A_{j_1}} \hat{d}_{i}\]
and
\[d_U=e+1>\frac{(k+1)e+e_s}{(k+1)d+d_s}d.\quad\]\qed

Fixing a coprime dimension vector $(d,e)$ we now deal with the question how to construct a certain set of stable $m$-bipartite quivers. Therefore, we assign a set of stable $m$-bipartite quivers to tuple of natural numbers which is uniquely determined by the dimension vector, see also Example $\ref{tupel}$. These numbers correspond to the number of possible glueing vertices and possible colourings of the constructed quivers.

Fix a dimension vector $(d,e)$ and the corresponding starting vector $(d_s,e_s)$. Denote by $\mathcal{T}^{(d,e)}_{n_1}$ the set of stable bipartite quivers of dimension type $(d_s,e_s)+n_1(d,e)$ with $n_1\geq 1$. As before let $\hat{\mathcal{T}}^{(d,e)}_{n_1}$ be the set which results by modifying a vertex $j_1$. Now we continue recursively: let $S\in\mathcal{T}^{(d,e)}_{n_k-1,\ldots,n_1}$ and $T\in\hat{\mathcal{T}}^{(d,e)}_{n_k,\ldots,n_1}$. Now let $\mathcal{T}^{(d,e)}_{1,n_k,\ldots,n_1}$ be the set consisting of all quivers $Q_{j_0,j_1}(S,T)$ such that $R_{j_0}+R_{j_1}\leq m$. Moreover, let the dimension of the glueing vertex $j$ be given by $\hat{d}_{j}=\hat{d}_{j_0}+\hat{d}_{j_1}-1$. In general let $\mathcal{T}^{(d,e)}_{n_{k+1},\ldots,n_1}$ 
be the set of glueing quivers resulting from glueing a quiver $S\in \mathcal{T}^{(d,e)}_{n_{k+1}-1,n_k,\ldots,n_1}$ and a quiver $T\in\hat{\mathcal{T}}^{(d,e)}_{n_k,\ldots,n_1}$ as described.
\begin{kor}\label{qst1}
The sets $\mathcal{T}^{(d,e)}_{n_k,\ldots,n_1}$ only contain stable quivers.
\end{kor}
{\it Proof.}
It suffices to prove that these quivers satisfy the conditions of Corollary $\ref{korverkl}$.

We assume that $\mathcal{T}^{(d,e)}_{n_k,\ldots,n_1}$ only contains stable quivers. We have to prove that $\mathcal{T}^{(d,e)}_{n_{k+1},\ldots,n_1}$ just consists of stable quivers for all $n_{k+1}\geq 1$. Therefore we show that the quivers in $\mathcal{T}^{(d,e)}_{n_k-1,\ldots,n_1}$ are starting quivers for all quivers in $\hat{\mathcal{T}}^{(d,e)}_{n_k,\ldots,n_1}$.

Let $(d^k,e^k)$ be the dimension type corresponding to $\hat{\mathcal{T}}^{(d,e)}_{n_k,\ldots,n_1}$ and $(d_s^k,e_s^k)$ the one belonging to $\mathcal{T}^{(d,e)}_{n_k-1,\ldots,n_1}$. It suffices to prove that
\[(d^{k+1}_s,e^{k+1}_s)=(d_s^k,e_s^k)+(n_k-1)(d^k,e^k)\]
is the starting vector for
\[(d^{k+1},e^{k+1})=(d_s^{k},e_s^{k})+n_k(d^{k},e^{k}).\]
Indeed, the quivers in $\hat{\mathcal{T}}^{(d,e)}_{n_k,\ldots,n_1}$ are obtained by the modification described in Corollary $\ref{korverkl}$. But this is equivalent to
\[e_s^{k+1}=\frac{1+d_s^{k+1}e^{k+1}}{d^{k+1}}\]
with the additional condition $d^{k+1}_s\leq d^{k+1}$, see Lemma $\ref{startexist}$. The second property follows immediately, the first one is equivalent to \[e_s^k=\frac{1+d_s^ke^k}{d^k},\]
what follows by a direct calculation. Therefore, the claim follows by the induction hypothesis.\qed
\begin{bei}\label{tupel}
\end{bei}
Let $(d_s,e_s)=(0,1)$ and $(d,e)=(1,n-1)$. Then we always obtain a corresponding tuple of natural numbers $(n_k,\ldots,n_1)$ to a fixed coprime dimension vector by proceeding as mentioned in Remark $\ref{darstellung}$. More detailed we have $(d^k,e^k)=(d_s,e_s)+n_k(d^{k-1},e^{k-1})$ and in this way we recursively obtain the whole tuple. The recursion terminates if $(d_s,e_s)=(0,1)$.

For instance consider $(d',e')=(5,8)$. The tuple of numbers is given by $(n_2,n_1)=(1,2)$ with $n=2$. Thus we get
\[(d', e') = (1, 2) + 2(2, 3) = (0, 1) + (1, 1) + (0, 1) + 2(1, 1)\]
Initially, consider the localization data of the dimension types $(1,2)$ and $(2,3)$, i.e.
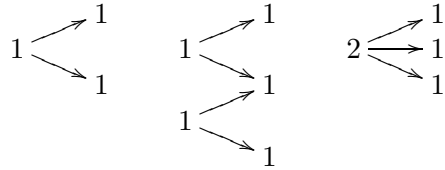
\[
\begin{xy}
\xymatrix@R0.5pt@C20pt{
&1&&1&&1\\1\ar[ru]\ar[rd]&&1\ar[ru]\ar[rd]&&2\ar[ru]\ar[r]\ar[rd]&1\\&1&&1&&1\\&&1\ar[ru]\ar[rd]&\\&&&1 }
\end{xy}
\]
where the numbers at the vertices indicate the dimension vector. By use of Corollary $\ref{korverkl}$ we obtain the following localization data of dimension type $(3,5)$ by glueing:
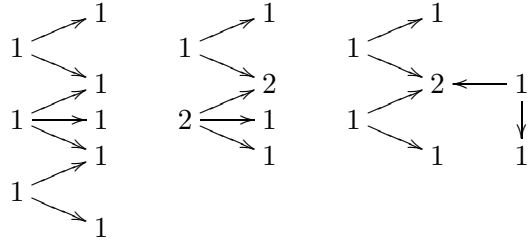
\[
\begin{xy}
\xymatrix@R0.5pt@C20pt{
&1&&1&&1\\1\ar[ru]\ar[rd]&&1\ar[ru]\ar[rd]&&1\ar[ru]\ar[rd]\\&1&&2&&2&1\ar[l]\ar[dd]\\1\ar[ru]\ar[rd]\ar[r]&1&2\ar[r]\ar[ru]\ar[rd]&1&1\ar[ru]\ar[rd]&\\&1&&1&&1&1\\1\ar[ru]\ar[rd]\\&1 }
\end{xy}
\]
Next, for instance we obtain the following localization data of type $(5,8)$ by glueing:
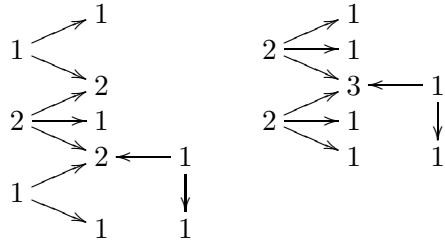
\[
\begin{xy}
\xymatrix@R0.5pt@C20pt{
&1&&&1\\1\ar[ru]\ar[rd]&&&2\ar[ru]\ar[r]\ar[rd]&1\\&2&&&3&1\ar[dd]\ar[l]\\2\ar[ru]\ar[rd]\ar[r]&1&&2\ar[ru]\ar[r]\ar[rd]&1\\&2&1\ar[dd]\ar[l]&&1&1\\1\ar[ru]\ar[rd]\\&1&1 }
\end{xy}
\]
\section{Asymptotics and combinatorics of trees}\label{strees}
\noindent The purpose of this section is to treat some aspects of combinatorics of trees. Fixing some properties we count the number of trees satisfying these properties, either exactly or at least asymptotically. This machinery will be used to count torus fixed points and fixed point components respectively. This gives rise to a lower bound for the number of fixed points and thus for the Euler characteristic of moduli spaces of the Kronecker quiver.\\

Let $a(x)=\sum_{n\geq 0}a_nx^n$ be a power series. In the following denote by
\[[x^n]a(x):=a_n\text{ with }n\geq 0\]
its $n$-th coefficient.
\begin{defi}
A tree is a connected acyclic graph. A rooted tree is a tree where a point is specified to be the root. A graph without cycles is called a forest, in particular, the components are trees.
\end{defi} 
When restricting to trees the points (resp. vertices) are often called knots. For further details according to trees and their combinatorics see for example \cite{hp} or \cite{sf}.
\subsection{Simply generated trees}\label{sgt}
\noindent We discuss simply generated trees, which we relate to localization data as constructed in the last section. Simply generated trees were introduced by Meir and Moon, see \cite{mm}, and are constructed as follows: fix a formal power series
\[\phi(x)=\sum_{n\geq 0}\phi_nx^n\]
such that $\phi_n\geq 0$ for all $n\geq 0$, $\phi_0>0$ and $\phi_j>0$ for at least one $j\geq 2$. Let $\mathcal{T}$ be the family of finite rooted trees. Define the weight $\omega_{\phi}(T)$ of a tree $T\in\mathcal{T}$ by
\[\omega_{\phi}(T)=\prod_{j\geq 0}\phi_j^{D_j(T)}\]
where $D_j(T)$ is the number of knots with $j$ successors. Denote by $|T|$ the number of knots of a tree $T$ and set
\[y_n=\sum_{|T|=n} \omega_{\phi}(T).\]
Now the generating function $y(x)=\sum_{n\geq 1}y_nx^n$ satisfies the functional equation
\[y(x)=x\phi(y(x)).\]Define $\mathcal{T}_{\phi}:=\{T\in\mathcal{T}\mid\omega_{\phi}(T)\neq 0\}$. We call a tree $T\in\mathcal{T}_{\phi}$ simply-generated by $\phi$.\\

For instance, if we define $\phi(x)=1+2x+x^2,$ we obtain the family of binary trees. Indeed, $y(x)$ satisfying $y(x)=x\phi(y(x))$ is its generating function, i.e. $y_n$ is the number of binary trees with $n$ knots. Here we take into account that we distinguish between left and right successors. 
\subsection{Lagrange inversion theorem}
\noindent In this section we briefly discuss the Lagrange inversion theorem, which will become an important tool later.
\begin{satz}\label{lagrange}
Let $\phi(x)=\sum_{n\geq 0}\phi_nx^n$ be a power series such that $\phi(0)\neq 0$ and let $y(x)$ be a power series satisfying the functional equation $y(x)=x\phi(y(x)).$
Let $g(x)$ be another power series. Then $y(x)$ is invertible and for the coefficients of $g(y(x))$ we have
\[[x^n]g(y(x))=\frac{1}{n}[u^{n-1}]g'(u)\phi(u)^n\]
for all $n\geq 1$.
Moreover, we have
\[[x^n](y(x))^m=\frac{m}{n}[u^{n-m}]\phi(u)^n.\]
\end{satz}

Note that this theorem is equivalent to the formulation of the Lagrange inversion theorem as usually stated in literature. For proofs and further details see for instance \cite{sf} or \cite{drm}.

By an easy calculation using Theorem $\ref{lagrange}$ we obtain the following special case which is important when counting localization data:
\begin{lem}
Let $\phi(x)=1+ax^b$ such that $y(x)=x\phi(y(x))$. Then we have 
\[[x^n]y(x)=\frac{1}{n}\binom{n}{\frac{n-1}{b}}a^{\frac{n-1}{b}}\]
if $b|n-1$ and $[x^n]y(x)=0$ otherwise.
\end{lem}
\begin{kor}
Let $m\geq 1$. We have
\[[x^n]y(x)^m=\frac{m}{n}\binom{n}{\frac{n-m}{b}}a^{\frac{n-m}{b}}\]
if $b|n-m$ and $n\geq m$ and $[x^n]y(x)^m=0$ otherwise.
\end{kor}

Let $a,b,m,n\in\mathbb{N}^+$.
Define
\[\mathcal{A}_{a,b,m,n}:=[x^n]y(x)^m\]
if $y(x)$ satisfies the functional equation $y(x)=x\phi(y(x))$ where $\phi(x)=1+ax^b$. Also define $\mathcal{A}_{a,b,n}:=\mathcal{A}_{a,b,1,n}$
\subsection{Asymptotic behaviour}
\noindent If we are not able to count the number of localization data exactly, we can use the following important result, see \cite{drm}, in order to count it asymptotically:
\begin{satz}\label{asym}
Let $F(x,y)$ be an analytic function in the variables $x$ and $y$ around $x=y=0$ such that $F(0,y)=0$ and such that the Taylor coefficients of $F$ around $0$ are real and non-negative. Then there exists an unique analytic solution $y=y(x)$ of the functional equation
\[y=F(x,y),\]
which has non-negative Taylor coefficients around $0$ and, moreover, $y(0)=0$.

If the region of convergence of $F(x,y)$ is large enough such that there exist positive solutions $x=x_0$ and $y=y_0$ of the system of functional equations given by \[y=F(x,y)\text{ and }1=F_y(x,y)\] with $F_x(x_0,y_0)\neq 0$ and $F_{yy}(x_0,y_0)\neq 0$, then $y(x)$ is analytic for $|x|<x_0$.

Moreover, there exist functions $h(x)$ and $g(x)$, which are  analytic around $x_0$, such that 
\[y(x)=g(x)-h(x)\sqrt{1-\frac{x}{x_0}}\]   
locally around $x_0$.

Then we have $g(x_0)=y(x_0)$ and
\[h(x_0)=\sqrt{\frac{2x_0F_x(x_0,y_0)}{F_{yy}(x_0,y_0)}}.\]
Furthermore, this provides a locally analytic continuation of $y(x)$ for $x-x_0\neq 0$.

If $[x^n]y(x)>0$ for all $n\geq n_0$, we also have that $x=x_0$ is the only singularity of $y(x)$ on the circle $|x|=x_0$. In conclusion, for $[x^n]y(x)$ we get an asymptotic expansion of the form
\[[x^n]y(x)=\sqrt{\frac{x_0F_x(x_0,y_0)}{2\pi F_{yy}(x_0,y_0)}}x_0^{-n}n^{-\frac{3}{2}}(1+\mathcal{O}(n^{-1})).\]
\end{satz}

Using the notations of the preceding theorem we get the following corollary:
\begin{kor}\label{x0}
Let $\phi(x)=1+ax^b$ and $y(x)$ such that $y(x)=x\phi(y(x))$. Then we have
\[(x_0)^{-1}:=ab\left( \frac{1}{(b-1)a}\right) ^{\frac{b-1}{b}}.\]
\end{kor}
\section{Applications}
\noindent In this section we discuss several consequences of the last sections and state several applications. First we discuss the asymptotic behaviour of the Euler characteristic of Kronecker moduli spaces. Then we consider some cases for which it is possible to calculate the Euler characteristic exactly.
\subsection{Conjecture concerning the asymptotic behaviour of the Euler characteristic}\label{verm}
\noindent In this subsection we discuss a conjecture based on Michael Douglas concerning the Euler characteristic of Kronecker moduli spaces and several consequences. Originally, in \cite{dou} Douglas suggested to fix $r\in\mathbb{R}_+$ and to consider $(d,e)\in\mathbb{N}^2_+$ with $\mathrm{gcd}(d,e)=1$ and $\frac{e}{d}\approx r$ to obtain the following:
\begin{enumerate}
\item There exists a $C_r\in\mathbb{R}$ such that for
$e,d\gg 0$ we have
\[\frac{\ln(\chi(M_{d,e}^m))}{d}\approx C_r.\] \item The function $r\mapsto
C_r$ is continuous. \end{enumerate}
Thus Douglas supposed that $\frac{\ln(\chi(M_{d,e}^m))}{d}$ and therefore the Euler characteristic is asymptotically already determined by the fraction $\frac{e}{d}$. Moreover, the Euler characteristic depends continuously on it. Let \[m_1:=\frac{m-\sqrt{m^2-4}}{2}\text{ and }m_2:=\frac{m+\sqrt{m^2-4}}{2}\,.\]
Based on this from \cite{wei} we obtain the following precise formulation:
\begin{ver}\label{ver}
There exists a continuous function
$f:[m_1,m_2]\subset\mathbb{R}\rightarrow\mathbb{R}$ such that the following holds: for all $r\in [m_1,m_2]$ and all $\varepsilon>0$ there exists an
$\delta>0$ and an $n\in\mathbb{N}$ such that for all
$(d,e)\in\mathbb{N}^2$ with $\gcd(d,e)=1$, $|r-e/d|<\delta$ and
$|d+e|>n$ we have
\[|f(r)-\frac{\ln(\chi(M^m_{d,e}))}{d}|<\varepsilon.\]
\end{ver}
\newpage\begin{bem}
\end{bem}
\begin{enumerate}
\item
We may also rephrase the conjecture as follows:
there exists a continuous function $f$ such that for every coprime dimension vector $(d,e)$ there exists a dimension vector $(d_s,e_s)$ such that
\[f(\frac{e}{d})=\lim_{n\rightarrow\infty}\frac{\ln{\chi(M_{d_s+nd,e_s+ne}^m)}}{d_s+nd}.\]
In particular, the right hand side converges.
\end{enumerate}

We discuss some consequences of the conjecture which are proved in \cite{wei}. For the remainder of this subsection we assume that the conjecture is true. Define
\[K:=(m-1)^2\ln((m-1)^2)-(m^2-2m)\ln(m^2-2m).\]
\begin{satz}
The function $f$ is given by
\[f(r)=\frac{K}{\sqrt{m-2}}\cdot\sqrt{r(m-r)-1}.\]
\end{satz}

In particular, the constant $K$ is its value at the point $r=1$. Moreover, we have that the Euler characteristic asymptotically only depends on the dimension of the moduli space:
\begin{kor}
The logarithm of the Euler characteristic $\ln(\chi(M_{d,e}^m))$ is asymptotically proportional to $\sqrt{dem-d^2-e^2}=\sqrt{\dim{M_{d,e}^m}-1}$.
\end{kor}

In this paper we prove that $f(1)=K$ and that the Euler characteristic grows exponentially. Note that, if a continuous function as conjectured exists, it follows from \cite{wei} that it is already uniquely determined by $f(1)$.
\subsection{The case of the dimension vector (d-1,d)}
\noindent In this section we investigate the function treated in Section \ref{verm} at the point $1$. This means investigating the dimension vector $(d-1,d)$. The Euler characteristic of the corresponding moduli space is, by applying the reflection functor, the same as the one corresponding to the dimension vector $(d,(m-1)d+1)$. The latter one is considered in the following. In particular, we show that the value at the point one is the one conjectured in Section \ref{verm}.\\

By Theorem \ref{bij} it is enough to consider the universal covering quiver of the Kronecker quiver $K(m)$. As a consequence, for the remainder of this subsection we only consider localization data such that the corresponding quiver is a subquiver of the universal covering quiver. Under this assumption, we will see that each localization data $(\mathcal{Q},\tilde{d})$ is of type one, i.e. $\tilde{d}_q\in\{0,1\}$ for all $q\in\mathcal{Q}_0$. This already implies that every localization data consists of subdata of dimension type $(1,m)$.
\begin{lem}
Every localization data $(\mathcal{Q},\tilde{d})$ of dimension type $(d,(m-1)d+1)$ is of type one. In particular, we have $\chi(M^s_{\tilde{d}}(\mathcal{Q}))=1$.
\end{lem}
{\it Proof.}
Let $(\mathcal{Q},\tilde{d})$ be a localization data of dimension type $(d,(m-1)d+1)$ and let $X$ be a stable representation of this data. Consider a subrepresentation
\[
\begin{xy}
\xymatrix@R1pt@C20pt{
&X_{j_{1}}\\&X_{j_{2}}\\X_i\ar[ruu]^{X_1}\ar[ru]_{X_2}\ar[rd]_{X_m}&\vdots\\&X_{j_{m}}\;}
\end{xy}
\]
The stability condition implies
\[d_{X_i}>\frac{(m-1)d+1}{d}\dim X_i>(m-1)\dim X_i .\]
In particular, this holds if $\dim X_i=\dim X_{j_k}=1$ for $k=1,\ldots,m$. Now we have $\dim X_{j_k} \geq \dim X_i $ for all $k$. Indeed, if we had $\dim X_{j_k} =l$ such that $l<\dim X_i $, we could consider the $(\dim X_i -l)$-subspace $\ker(X_k)$ which would just have a $(\dim X_i -l)(m-1)$-dimensional image. This contradicts the stability condition.

Therefore, the subrepresentation is of dimension type $(\dim X_i ,e')$ with $e'\geq m\dim X_i$. 
Furthermore, because of the stability every $k$-dimensional subspace at least has an $((m-1)k+1)$-dimensional image.

If we fix a proper boundary quiver, which exists because the original quiver has no cycles, this subquiver just has one common vertex with the remainder of the quiver and the corresponding subdata is of dimension type $(d_1,md_1)$. But for the dimension type $(d-d_1,b)$ of the remainder of the data we have
\[b\geq (m-1)(d-d_1)+1.\]  
Let $h\geq 1$ be the dimension of the intersection of (the vector spaces corresponding to the common vertex of) the two subrepresentations of $X$ corresponding to the two subdata. Then we get
\[(m-1)d+1=b+d_1m-h\geq (m-1)(d-d_1)+1+d_1m-h=(m-1)d+d_1-h+1.\]
Therefore, we have $h\geq d_1$ and thus $h=d_1$.

We continue by proving that after removing the subdata of dimension type $(d_1,(m-1)d_1)$, i.e. the subdata of dimension type $(d_1,md_1)$ except the common vertex, we get a localization data of dimension type $(d-d_1,(m-1)(d-d_1)+1)$. It suffices to prove stability because the original subdata has a $d_1$-dimensional intersection with the remainder.

For an arbitrary subspace $U\subset \oplus_{i\in I}X_i$ with $\dim U<d-d_1$ we have
\[d_U>\frac{(m-1)d+1}{d} \dim U.\]
Since $\dim U<d-d_1$, we also have
\[d_U>\frac{(m-1)(d-d_1)+1}{(d-d_1)} \dim U\]
proving the claim in-between.

Thus we can proceed by induction on the number of sources in order to show that all localization data are of type one.

Consider some data such that corresponding quiver has one source. Obviously, it is a stable quiver of type $1$.

Assume that the quiver has $n+1$ sources. We may remove a proper boundary quiver so that we again get a localization data, which is of the requested type by induction hypothesis. But since the original quiver has no cycles, there exist at least two proper boundary quivers. Thus the assertion follows by applying the induction hypothesis to the respective subquivers after removing a proper boundary quiver.
The second statement for instance follows when considering the dimension formula mentioned in Remark \ref{bem1}.\qed
\begin{satz}
We have
\[\chi(M_{d,d+1}^m)=\frac{m}{(d+1)((m-1)d+m)}\binom{(m-1)^2d+(m-1)m}{ d}.\]
Moreover, we also have
\[f(1)=\lim_{d\rightarrow\infty}\frac{\ln(\chi(M_{d,d+1}^m)}{d}=(m-1)^2\ln(m-1)^2-(m^2-2m)\ln(m^2-2m)\]
for $f$ defined in Section \ref{verm}.
\end{satz}
{\it Proof.}
As shown previously, we may assume that all subdata of a localization data with one source have vertex set \[I\cup J=\{i,j_1,\ldots,j_m\}\] and arrow set \[R=\{(i,j_1),\ldots,(i,j_m)\}\]
with $\tilde{d}_{i}=\tilde{d}_{j_k}=1$. In particular, the moduli spaces of the considered localization data are zero-dimensional yielding that the Euler characteristic is one.

By Remark $\ref{anzko}$ there exists exactly one possibility to choose a colouring $c$ taking into account the symmetries of $S_m$. Again by Remark $\ref{anzko}$ and by Corollary \ref{korverkl} we can glue $k$ subquivers on each vertex $j_l$, $1\leq l\leq m$, with $0\leq k\leq (m-1)$ in order to get a localization data. But we have to take note of the symmetries of $S_k$. Assuming that there is only one starting knot let $y(x)$ the generating function of such quivers and consider
\begin{eqnarray*}\phi(x)&=&1+\frac{(m-1)}{|S_1|}x^{m-1}+\frac{(m-1)(m-2)}{|S_2|}x^{2(m-1)}\ldots+\frac{\prod_{i=1}^{m-1}(m-i)}{|S_{m-1}|}x^{(m-1)(m-1)}\\&=&\sum_{i=0}^{m-1}x^{i(m-1)}\binom{m-1}{i}=(1+x^{m-1})^{m-1}.\end{eqnarray*}
By Section $\ref{sgt}$ the generating function satisfies the functional equation $y(x)=x(\phi(y(x)))$. Now the generating function for all localization data is obtained as follows: we start with the unique localization data of dimension type $(1,m)$ having $m$ knots. The resulting generating function is $y(x)^m$ and by applying the Lagrange inversion theorem we obtain that
\[[x^n]y(x)^m=\frac{m}{n}[u^{n-m}]\phi(u)^n=\frac{m}{n}\binom{n(m-1)}{\frac{n-m}{m-1}}.\]
If we assign the weight $0$ to the sink of the starting quiver, every such quiver that has $(m-1)d+1$ knots corresponds to a localization data of dimension type $(d,(m-1)d+1)$. The other way around, we may assume that every localization data has some sink $i$ with weight $0$ what gives us $d$ choices. This means for every localization data we exactly get $d$ trees. Hence we get
\begin{eqnarray*}\chi(M_{d,(m-1)d+1}^m)&=&\frac{m}{d((m-1)d+1)}\binom{(m-1)^2d+(m-1)}{ d-1}\\
&=&\frac{m}{d((m-1)(d-1)+m)}\binom{(m-1)^2(d-1)+(m-1)m}{ d-1}.\\
\end{eqnarray*}
Since $\chi(M_{d-1,d}^m)=\chi(M_{d,(m-1)d+1}^m)$, the assertion is proved.\\

The second part follows by applying Theorem $\ref{asym}$. It may be left unconsidered that exactly $d$ trees define the same localization data. Indeed, obviously we have
\[\lim_{d\rightarrow\infty}\frac{\ln{d}}{d}=0.\]
We can also assume that we just have one starting knot. Thus in addition to the functional equation $y(x)=x(\phi(y(x)))$ we consider the functional equation
\[1=x(m-1)^2(1+y(x)^{m-1})^{m-2}y(x)^{m-2}.\]
Moreover, we consider the equations
\[y_0=x_0(1+y_0^{m-1})^{m-1}\]
and
\[1=x_0(m-1)^2(1+y_0^{m-1})^{m-2}y_0^{m-2}.\]
Then we have
\[x_0=\frac{1}{(m-1)^2(1+y_0^{m-1})^{m-2}y_0^{m-2}}\]
which implies
\[y_0=\frac{(1+y_0^{m-1})^{m-1}}{(m-1)^2(1+y_0^{m-1})^{m-2}y_0^{m-2}}.\]
Thus we get
\[(m-1)^2y_0^{m-1}=1+y_0^{m-1}\]
and therefore
\[y_0^{m-1}=\frac{1}{(m-1)^2-1}.\]
Hence we have
\begin{eqnarray*}(x_0)^{-1}&=&(m-1)^2\left(\frac{(m-1)^2}{(m-1)^2-1}\right)^{m-2}\left(\frac{1}{(m-1)^2-1}\right)^{\frac{m-2}{m-1}}\\
&=&(m-1)^{2(m-1)}\left(\frac{1}{m^2-2m}\right)^{m-2+\frac{m-2}{m-1}}=(m-1)^{2(m-1)}\left(\frac{1}{m^2-2m}\right)^{\frac{m^2-2m}{m-1}}.\\
\end{eqnarray*}
If the number of knots of the considered trees is $(m-1)d+1$, such a tree corresponds to a localization data of dimension type $(d,(m-1)d+1)$ for all $d\geq 1$.
Since we consider the logarithm, we may discount the remaining factors of Theorem $\ref{asym}$. Thus we obtain
\begin{eqnarray*}\lim_{d\rightarrow\infty}\frac{\ln\chi(M_{d,(m-1)d+1}^m)}{d}&=&\lim_{d\rightarrow\infty}\frac{\ln\left((m-1)^{2(m-1)}\left(\frac{1}{m^2-2m}\right)^{\frac{(m^2-2m)}{m-1}}\right)^{(m-1)d+1}}{d}\\
&=&(m-1)^2\ln(m-1)^2-(m^2-2m)\ln(m^2-2m).
\end{eqnarray*}
Because of the isomorphisms of moduli spaces, the assertion follows.\qed

\subsection{A lower bound}\label{62}
\noindent The aim of this section is to determine a lower bound for the Euler characteristic of Kronecker moduli spaces for coprime dimension vectors which also proves the exponential growth of the Euler characteristic as conjectured by Douglas. Therefore, we consider such dimension vectors $(d,e)$ of $K(m)$ satisfying $e>(m-1)d$. The remaining cases are obtained by the isomorphisms of the moduli spaces stated in Proposition $\ref{propofkron}$. In the considered cases the moduli spaces are zero-dimensional. Moreover, we will see that the recursive construction of the localization data simplifies.\\

As a consequence of Theorem $\ref{bij}$ we again assume that all torus fixed points are representations of the universal covering quiver.

Initially, consider the dimension vectors $(1,n-1)$ and $(1,n)$ with $2\leq n\leq m-1$ which correspond to the dimension vectors $(n-1,m(n-1)-1)$ and $(n,mn-1)$ by the mentioned isomorphisms. For the dimension type $(1,n-1)$ there exists only one localization data
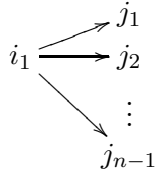
\[
\begin{xy}
\xymatrix@R0.5pt@C20pt{
&j_1\\i_1\ar[ru]\ar[r]\ar[rdd]&j_2\\&\vdots\\&j_{n-1}}
\end{xy}
\]
where $\tilde{d}_{j_k}=\tilde{d}_{i_1}=1$ for all $1\leq k\leq n-1$.
Analogously, we obtain the unique localization data of dimension type $(1,n)$. 

Consider the following localization data of dimension type $(n-1,m(n-1)-1)$ where $\tilde{d}_j=n-2$ and $\tilde{d}_{j_{k,l}}=\tilde{d}_{i_{k}}=1$ otherwise:
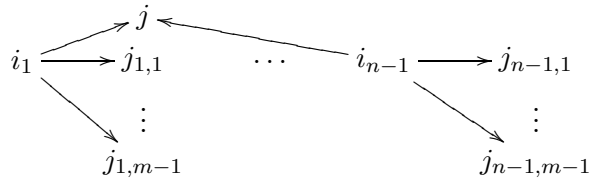
\[
\begin{xy}
\xymatrix@R0.5pt@C20pt{
&j&\\i_1\ar[ru]\ar[r]\ar[rdd]&j_{1,1}&\dots&i_{n-1}\ar[rdd]\ar[llu]\ar[r]&j_{n-1,1}\\&\vdots&&&\vdots\\&j_{1,m-1}&&&j_{n-1,m-1}}
\end{xy}
\]
Again we analogously obtain the data of type $(n,mn-1)$.
\begin{bem}
\end{bem}
\begin{enumerate}
\item For the dimension vector $(n,mn-1)$, $1\leq n\leq m$, this is also the only localization data because obviously each one-dimensional subspace is forced to have an $m$-dimensional image. Moreover, because of the stability condition, we have for each subspace $U$ of dimension $d'<n$ which corresponds to a stable representation of this data that
\[d_U>\frac{nm-1}{n}d'.\]
Therefore, we have $d_U\geq md'$ for all $d'<n$. But, for any other data of this dimension type this condition is not satisfied.
\item We also get this localization data by applying the reflection functor, see Theorem $\ref{kspi}$.
\end{enumerate}

By use of the procedure introduced in Section $\ref{stabilitaet}$ we can glue these quivers. Fix $m\in\mathbb{N}$ and define $Q^l$ by
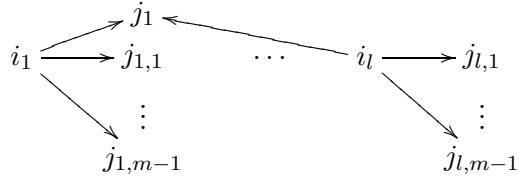
\[
\begin{xy}
\xymatrix@R0.5pt@C20pt{
&j_1&\\i_1\ar[ru]\ar[r]\ar[rdd]&j_{1,1}&\dots&i_{l}\ar[rdd]\ar[llu]\ar[r]&j_{l,1}\\&\vdots&&&\vdots\\&j_{1,m-1}&&&j_{l,m-1}}
\end{xy}
\]
Let $I\cup J$ be the set of vertices and define $J_1':=J\backslash j_1$. Let $\tilde{d}_{j_1}=l-1$ and let $\tilde{d}_q=1$ for the remaining vertices. Define the glueing quiver $Q^{l_1,l_2}:=Q_{j,j_2}(Q^{l_1},Q^{l_2})$ with $j\in J'_1$. For the resulting data define $\tilde{d}_{j_2}=l_2$ whereby the dimensions of the other vertices remain constant. For instance we obtain:
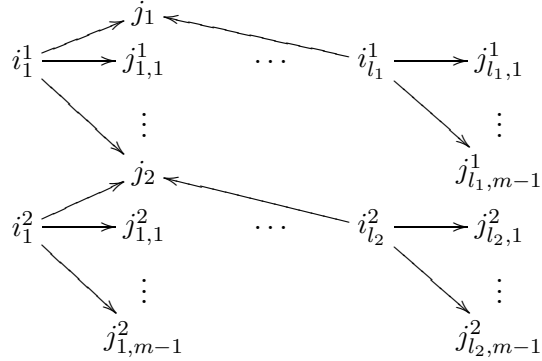
\[
\begin{xy}
\xymatrix@R0.5pt@C20pt{&j_1&\\i_1^1\ar[ru]\ar[r]\ar[rdd]&j^1_{1,1}&\dots&i^1_{l_1}\ar[rdd]\ar[llu]\ar[r]&j^1_{l_1,1}\\&\vdots&&&\vdots\\&j_2&&&j^1_{l_1,m-1}\\i^2_1\ar[ru]\ar[r]\ar[rdd]&j^2_{1,1}&\dots&i^2_{l_2}\ar[rdd]\ar[llu]\ar[r]&j^2_{l_2,1}\\&\vdots&&&\vdots\\&j^2_{1,m-1}&&&j^2_{l_2,m-1}}
\end{xy}
\]
We again consider the construction of Corollary $\ref{qst1}$. Let $(d,e)=(n_1(n-1)+1,n_1(m(n-1)-1)+m)=(1,m)+n_1(n-1,m(n-1)-1)$ with $n_1\in\mathbb{N}$. Then we obtain the cases \begin{equation}\label{eig}\frac{m(n-1)-1}{n-1}d\leq e\leq\frac{mn-1}{n}d.\end{equation} Now the quivers are glued as explained in Section \ref{stabilitaet}. Fixing $n_1\geq 1$ we denote the resulting data by $\mathcal{Q}^{n}_{n_1}$. They obviously result if one successively glues $n_1$-times some data of dimension type $(n-1,m(n-1)-1)$ to some data of type $(1,m)$. Call the glueing vertex corresponding to the first glueing initial glueing vertex. If $j_1$ is the initial glueing vertex, denote by $\hat{\mathcal{Q}}^{n}_{n_1}$ the set of quivers obtained by increasing the dimension of $j_1$ by one.

We now recursively define
\[\mathcal{Q}^n_{n_{k+1},\ldots,n_1}=\{Q_{j,j_1}(S,T)\mid S\in\mathcal{Q}^n_{n_{k+1}-1,n_k,\ldots,n_1},T\in\hat{\mathcal{Q}}^n_{n_k,\ldots,n_1}\},\]
where $j\in S_0$ such that $R_j=1$ and where $j_1$ is the initial glueing vertex of $T\in\mathcal{Q}^n_{n_k,\ldots,n_1}$. Furthermore, let $\mathcal{Q}^n_{0,n_k,\ldots,n_1}=\mathcal{Q}^n_{n_k-1,\ldots,n_1}$. By Corollary $\ref{qst1}$ we know that every data $S\in\mathcal{Q}^n_{n_{k+1}-1,\ldots,n_1}$ is a localization data and that each of them satisfies the properties of the starting quiver for each $T\in\hat{\mathcal{Q}}^n_{n_k,\ldots,n_1}$. Thus it follows that every data which is obtained in such a way is a localization data.
\begin{bem}
\end{bem}
\begin{enumerate}
\item If $(d,e)$ is given such that $(\ref{eig})$ holds, we can determine the corresponding tuple $(n_k,\ldots,n_1)$ as described in Remark $\ref{darstellung}$. Note that there is an easier method to get this tuple by simply solving linear equations, see \cite{wei3}.
\end{enumerate}

Next we determine the cardinality of these sets in order to obtain a lower bound for the Euler characteristic. The moduli spaces of the considered localization data are zero-dimensional, i.e. a point. Furthermore, by Theorem $\ref{propofkron}$ we can assume that $n\geq\frac{m+1}{2}$. This is another advantage simplifying combinatorics. Indeed, because of this assumption it is just possible to glue {\it one} quiver on each vertex of dimension one. Otherwise, there is no suitable colouring to obtain a localization data from the produced quiver because it is no subquiver of the regular $m$-tree.\\

Initially, consider the set $\mathcal{Q}_1$ consisting of the localization data of dimension type $(n,mn-1)$. After modifying a sink, considering the properties of Remark $\ref{anzko}$ and taking into account all symmetries and the fact that all quivers are glued as mentioned above, there exist 
\[\binom{m-1}{ n}\text{ possibilities}\]  
to choose a colouring $c:R\mapsto \{1,\ldots,m\}$ where $R$ is the set of arrows.

Each of the quivers has $n(m-1)$ knots, i.e. vertices $j\in J$ such that $R_j=1$. Denote by  $a^n_{n_1}$ the cardinality of $\hat{\mathcal{Q}}^n_{n_1}$ in consideration of the different colourings. Furthermore, let $K^n_{n_1}$ the number of knots of these quivers which coincide for all quivers in this set.

Using the notation of Section $\ref{strees}$ we have 
\[a^n_{n_1}=\binom{m-1}{n}\mathcal{A}_{\binom{m-1}{ n-1},(n-1)(m-1),n(m-1),n(m-1)+(n_1-1)(n-1)(m-1)}.\]
Moreover, we have
\[K^n_{n_1}=n(m-1)+(n_1-1)(n-1)(m-1)-(n_1-1).\]
Considering the construction we get the following lemma by an easy observation.
\begin{lem}\label{kno}
Let $(n_{k+1},\ldots,n_1)\in\mathbb{N}^{k+1}$.
\begin{enumerate}
\item The number of knots of the quivers in $\hat{\mathcal{Q}}^n_{n_{k+1},\ldots,n_1}$ is given by
\[K^n_{n_{k+1},\ldots,n_1}=K^n_{n_{k}-1,\ldots,n_1}+n_{k+1}K^n_{n_{k},\ldots,n_1}-n_{k+1}.\]
\item 
Moreover, we have
\[a^n_{n_{k+1},\ldots,n_1}=a^n_{n_{k}-1,\ldots,n_1}\cdot\mathcal{A}_{a^n_{n_{k},\ldots,n_1},K^n_{n_{k},\ldots,n_1},K^n_{n_{k}-1,\ldots,n_1},
K^n_{n_{k}-1,\ldots,n_1}+n_{k+1}K^n_{n_{k},\ldots,n_1}}.\]
\end{enumerate}
\end{lem}

Fixing a dimension vector, it suffices to determine the corresponding tuple of natural numbers in order to get a lower bound for the Euler characteristic. Given a tuple as above define $K_{d,e}^m:=K^n_{n_{k+1},\ldots,n_1}$ and  $a_{d,e}^m:=a^n_{n_{k+1},\ldots,n_1}$ and consider the function \[\phi(x)=1+a_{d,e}^mx^{K_{d,e}^m}.\]
The generating function $y(x)$ satisfies the functional equation $y(x)=x\phi(y(x))$.
Since we are interested in some asymptotic value, which is independent of the number of starting knots, we can assume that there exists just one starting knot. Even the starting quiver only gives us a constant, which we may ignore.\\

For every coloured tree constructed like this we obtain some localization data by assigning the weight $0$ to the source of the starting quiver. Thus it may happen that different trees define the same localization data. But, if $(d,e)$ is the considered dimension vector, the number of possible starting quivers is bounded by $d$. Since
\[\lim_{d\rightarrow\infty}\frac{\ln{d}}{d}=0,\] we may disregard this as well when investigating the logarithmic asymptotic behaviour. 
Define \[u_{d,e}^m:=\frac{K_{d,e}^m}{d}.\]
\begin{satz}
Let $e>(m-1)d$. We have
\[\lim_{n\rightarrow\infty}\frac{\ln(\chi(M_{d_s+nd,e_s+nd}))}{d_s+nd}\geq \frac{1}{d}(\ln a_{d,e}^m+K^m_{d,e}\ln K^m_{d,e}-(K^m_{d,e}-1)\ln (K^m_{d,e}-1)).\]
\end{satz}
{\it Proof.}
Define 
\[F(x,y)=x\phi(y(x)).\]
By Corollary $\ref{x0}$ there exists a constant $C\in\mathbb{R}_{>0}$ such that $[x^n]y(x)=Cx_0^{-n}n^{-\frac{3}{2}}(1+\mathcal{O}(n^{-1}))$
with \[(x_0)^{-1}=a_{d,e}^mK_{d,e}^m\left( \frac{1}{(K_{d,e}^m-1)a_{d,e}^m}\right) ^{\frac{K_{d,e}^m-1}{K_{d,e}^m}}.\]
Then we get
\begin{eqnarray*}((x_0)^{-1})^{nK_{d,e}}&=&(a_{d,e}^mK_{d,e}^m)^{nu_{d,e}^md}\left( \frac{1}{(K_{d,e}^m-1)a_{d,e}^m}\right)^{nu_{d,e}^md-n}\\&&=(a_{d,e}^m(K_{d,e}^m-1))^n\left(\frac{K_{d,e}^m}{(K_{d,e}^m-1)}\right)^{nu_{d,e}^md}.
\end{eqnarray*}
Hence we get that
\begin{align*}&\frac{\ln(\chi(M_{d_s+nd,e_s+nd}))}{d_s+nd}&\\&\geq \ln K+\frac{n\cdot \ln(a_{d,e}^m(K_{d,e}^m-1))}{d_s+nd}+\frac{nu_{d,e}^md \cdot \ln(K_{d,e}^m)}{d_s+nd}-\frac{nu_{d,e}^md\cdot \ln(K_{d,e}^m-1)}{d_s+nd}=:L_{d,e,n}^m
\end{align*}
for a constant $K\in\mathbb{R}_{>0}$. Thus it follows
\begin{eqnarray*}\lim_{n\rightarrow\infty}L_{d,e,n}^m&=&\frac{\ln(a_{d,e}^m(K_{d,e}^m-1))}{d}
+u_{d,e}^m(\ln K_{d,e}^m-\ln(K_{d,e}^m-1)\\
&=&\frac{1}{d}(\ln a_{d,e}^m+K^m_{d,e}\ln K^m_{d,e}-(K^m_{d,e}-1)\ln (K^m_{d,e}-1))=:L_{d,e}^m
\end{eqnarray*}
which proves the theorem.\qed

By use of the isomorphisms of the moduli spaces we also get a lower bound for arbitrary $d$ and $e$.
\begin{bei}
\end{bei}
This example applies the introduced methods to the case $(d,e)=(5,8)$ and $m=3$. For the starting dimension vector we get $(d_s,e_s)=(3,5)$, for the localization data of this type see Example \ref{tupel}.

The reflected dimension vector is $(8,19)$ and we obtain $K_{5,8}^3=12$ and $a_{5,8}^3=1664$. Thus in conclusion we have 
\[L_{5,8}^3=\frac{1}{5}\ln{\left(1664\cdot\frac{12^{12}}{11^{11}}\right)}.\]
\subsection{The case of the dimension vector (3,4)}\label{34}
\noindent In this section we consider the case $d=3$ and $e=4$ with $m\geq 3$ in detail. Consider the stable bipartite quiver given by
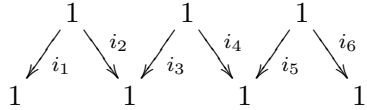
\[
\begin{xy}
\xymatrix@R2pt@C10pt{
&1\ar[lddd]^{i_1}\ar[rddd]^{i_2}&&1\ar[lddd]^{i_3}\ar[rddd]^{i_4}&&1\ar[lddd]^{i_5}\ar[rddd]^{i_6}&\\\\\\1&&1&&1&&1}
\end{xy}
\]
Therefore, by colouring the arrows in the colours $\{1,\ldots,m\}$ satisfying the conditions of Remark $\ref{anzko}$ we obtain a localization data. In this case, the conditions are $c(i_l)\neq c(i_{l+1})$ for $1\leq l\leq 5$. Each colouring is unique up to the symmetry of the symmetric group $S_2$.

The colourings $(i,j,k,i,j,k)$ and $(i,j,k,i,j,i)$, such that $i,j,k\in\{1,\ldots,m\}$ are pairwise disjoint, give rise to two cases, which we now consider in greater detail. In the first case we obtain
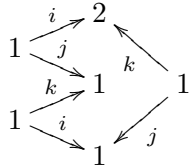
\[
\begin{xy}
\xymatrix@R0.5pt@C20pt{
&2&\\1\ar[ru]^{i}\ar[rd]^{j}&&\\&1&1\ar[ldd]^{j}\ar[luu]^k\\1\ar[ru]^{k}\ar[rd]^{i}&&\\&1&}
\end{xy}
\]
There is no new symmetry arising from this colouring. Furthermore, the moduli space is a point for this dimension vector. Note that the cycle breaks down after a second localization so that we get back the former quiver.

The second special case is
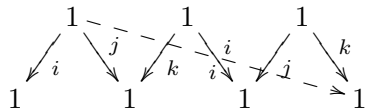
\[
\begin{xy}
\xymatrix@R2pt@C10pt{
&1\ar@{-->}[rrrrrddd]_i\ar[lddd]^{i}\ar[rddd]^{j}&&1\ar[lddd]^{k}\ar[rddd]^{i}&&1\ar[lddd]^{j}\ar[rddd]^{k}&\\\\\\1&&1&&1&&1}
\end{xy}
\]
The colouring induces an extra arrow and therefore another symmetry. In particular, the localization data is already determined by the choice of the colour of the free arrow, i.e. the one that does not appear in the cycle. But because of the extra arrow the moduli space is  $\mathbb{P}^1$ so that the Euler characteristic is two.

Note that $\chi(\mathbb{P}^1)=2$ follows also from a second localization. Indeed, by considering the quiver without its colouring the fixed points are those representations satisfying $X_{i5}=0$ or $X_{i7}=0$ where $i_7$ is the extra arrow. Thus we again get back the original localization data by a second localization. In conclusion we obtain that there are $\frac{m(m-1)^5}{|S_2|}$ possibilities to choose a colouring.

Further localization data are given by colourings of the following stable bipartite quiver:
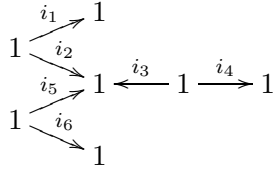
\[
\begin{xy}
\xymatrix@R0.5pt@C20pt{
&1\\1\ar[ru]^{i_1}\ar[rd]^{i_2}&\\&1&1\ar[l]_{i_3}\ar[r]^{i_4}&1\\1\ar[ru]^{i_5}\ar[rd]^{i_6}&\\&1 }
\end{xy}
\]
with the conditions $c(i_1)\neq c(i_2)$, $c(i_3)\neq c(i_4)$, $c(i_5)\neq c(i_6)$ and $c(i_2)$, $c(i_3)$, $c(i_5)$ pairwise disjoint. In consideration of the symmetries of $S_3$ we obtain $\frac{m(m-1)^4(m-2)}{|S_3|}$ possibilities. We also get
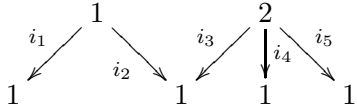
\[
\begin{xy}
\xymatrix@R2pt@C20pt{
&1\ar[lddd]_{i_1}\ar[rddd]_{i_2}&&2\ar[lddd]_{i_3}\ar[ddd]^{i_4}\ar[rddd]^{i_5}&\\\\\\1&&1&1&1}
\end{xy}
\]
with the conditions $c(i_3)$, $c(i_4)$, $c(i_5)$ pairwise disjoint and $c(i_1)\neq c(i_2)\neq c(i_3)$. Thus we get $\frac{m(m-1)^3(m-2)}{|S_2|}$ possibilities.\\

If $m\geq 4$, we finally get the localization data coming from
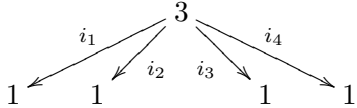
\[\begin{xy}
\xymatrix@R2pt@C20pt{
&&3\ar[llddd]_{i_1}\ar[lddd]^{i_2}\ar[rddd]_{i_3}\ar[rrddd]^{i_4}&&\\\\\\1&1&&1&1}
\end{xy}
\]with the condition that the colours of all arrows are pairwise disjoint, hence $\binom{m}{ 4}$ possibilities.

Since all fixed point components may be understood as points, for the Euler characteristic we have \[\chi(M_{3,4}^m)=\binom{m}{4}+\frac{m(m-1)^3(m-2)}{2}+\frac{m(m-1)^4(m-2)}{6}+\frac{m(m-1)^5}{2}.\]
One easily verifies that this is the same result one obtains by the algorithm from \cite{rei2}, i.e.:
\[\chi(M_{3,4}^m)=\frac{1}{24}m(m-1)(4m^2-7m+2)(4m^2-7m+1).\]

\subsection{The case of the dimension vector $(d,d)$}
\noindent The next application is to consider the Euler characteristic of Kronecker moduli spaces corresponding to the dimension vectors $(d,d)$, $d\in\mathbb{N}$. We will see that the Euler characteristic vanishes if $d\geq 2$. In this section we consider the Kronecker quiver $K(m)$ with $m\geq 1$.
\begin{lem}
Every stable torus fixed point $X=((V,W),(X_1,\ldots,X_m))$ of $M^m_{d,(m-1)d}$ has a cycle. Thus there exists a subspace $U\subset W$ and maps $f_1,\ldots,f_{2k}\in\{X_1,\ldots,X_m\}$ with $f_i\neq f_{i+1}$ for $1\leq i\leq 2l-1$ such that
\[f_1\circ f_2^{-1}\ldots\circ f_{2k-1}\circ f_{2k}^{-1}(U)=U.\]
\end{lem}
\begin{bem}
\end{bem}
\begin{enumerate}
\item
From the proof we even get the stronger result that the quiver of some localization data with this dimension is forced to be cyclic. In particular, there exists no subquiver having just one common vertex with the remainder of the quiver.
\end{enumerate}
{\it Proof.}
Let $(\mathcal{Q},\hat{d})$ be a localization data and let $X$ be a stable representation of this data. Consider a subdata of the form
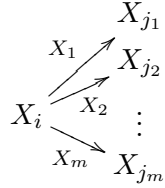
\[
\begin{xy}
\xymatrix@R1pt@C20pt{
&X_{j_{1}}\\&X_{j_{2}}\\X_i\ar[ruu]^{X_1}\ar[ru]_{X_2}\ar[rd]_{X_m}&\vdots\\&X_{j_{m}}}
\end{xy}
\]
Because of the stability we have
\[d_{X_i}>\frac{(m-1)d}{d}\dim X_i =(m-1)\dim X_i.\]
We also have $\dim X_{j_k} \geq \dim X_i $ for all $k$. Indeed, if we had $\dim X_{j_k} =l$ such that $l<\dim X_i $, we could consider the $(\dim X_i-l)$-subspace $\ker(X_k)$. It would have a $(\dim X_i -l)(m-1)$-dimensional image, which obviously contradicts the stability condition.

Therefore, the subdata is of dimension type $(\dim X_i ,e')$ with $e'\geq m\dim X_i $. 
Moreover, the stability implies that each $k$-dimensional subspace corresponding to $X$ has at least an $((m-1)k+1)$-dimensional image.

Assume that the localization data would not have a cycle. Thus, in particular, it would have some proper boundary quiver which apparently would be of dimension type $(d_1,md_1)$. If we denote by $b$ the sum of the dimensions corresponding to the sinks of the remainder of the data, we get
\[b\geq (m-1)(d-d_1)+1.\]  
Define $h:=(m-1)(d-d_1)+b-(m-1)d$ which is the minimal possible dimension of the intersection of the two corresponding subrepresentations (at the common vertex) of a stable representation of the considered data. It follows
\[(m-1)d=b+d_1m-h\geq (m-1)(d-d_1)+1+d_1m-h=(m-1)d+d_1-h+1.\]
It follows $h\geq d_1+1$ and thus $d_1=0$.\qed
\begin{kor}\label{vanish}
The Euler characteristic of the Kronecker moduli spaces with dimension vector $(d,d)$ vanishes if $d\geq 2$.
\end{kor}
{\it Proof.}
By the previous lemma we know that each representation of a localization data of dimension type $(d,(m-1)d)$ has a cycle. But because of Theorem $\ref{bij}$ we can assume that fixed points of each Kronecker moduli space do not have cycles. Hence there are no stable representations of the universal covering quiver of dimension type $(d,(m-1)d)$.

Because of the isomorphism between $M_{d,(m-1)d}^m$ and $M_{d,d}^m$, we in conclusion get
\[\chi(M_{d,(m-1)d}^m)=\chi(M_{d,d}^m)=0\quad.\]\qed
\subsection{Finiteness of the fixed point set}
\noindent In this section we investigate and answer a question posed in \cite{dre}. Namely for which coprime dimension vectors is the set of fixed points finite and for which dimension vectors exists at least one $n$-dimensional fixed point component with $n\geq 1$. 
\begin{satz}
Let $d\geq 3$, $e\geq 4$ and $m\geq 3$ . Then there exist infinitely many torus fixed points.
\end{satz}
{\it Proof.} Since the torus action is compatible with the isomorphisms, we may assume
\[d\leq e\leq\frac{m}{2}d.\]
Furthermore, let $m'\leq \frac{m}{2}\in\mathbb{N}$ such that $(m'-1)d<e<m'd$. By \cite{wei3} there exists a stable bipartite quiver $s_{d,e}^m$ of type one which consists of subdata of dimension type $(1,m')$ and $(1,m'+1)$ respectively. Since $d\geq 3$, there exists a subdata of the form
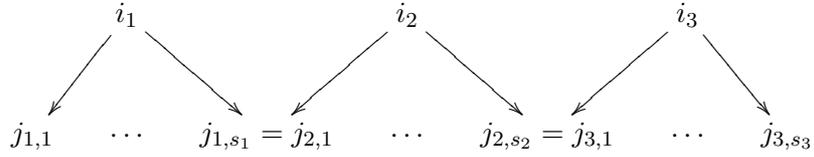
\[
\begin{xy}
\xymatrix@R6pt@C15pt{
&i_1\ar[lddd]\ar[rddd]&&i_2\ar[lddd]\ar[rddd]&&i_3\ar[lddd]\ar[rddd]&\\\\\\j_{1,1}&\cdots&j_{1,s_1}=j_{2,1}&\cdots&j_{2,s_2}=j_{3,1}&\cdots&j_{3,s_3}}
\end{xy}
\]
with $s_1,s_2,s_3\in\{m',m'+1\}$. Fix an arbitrary colouring $c$ of the arrows which satisfies
\[c(i_1,j_{1,1})=c(i_2,j_{2,s_2})=c(i_3,j_{3,s_3})=1, c(i_1,j_{1,s_1})=c(i_3,j_{3,1})=2\text{ and }c(i_2,j_{2,1})=3\]
and $c(i_3,j_{3,k})\neq 3$ for every $k=1,\ldots, s_3$. This is possible because $s_3<m$. This colouring induces an extra arrow $(i_3,j_{1,1})$ such that $c(i_3,j_{1,1})=3$. Hence the associated moduli space is at least one-dimensional implying that there are infinitely many torus fixed points.\qed

\subsection{Open questions}
\noindent A fundamental question is how to determine all localization data and if it is perhaps enough to know all localization data of type one. Also, one could ask if it is possible to put all localization data down to the case of localization data of type one. For instance, when considering the stable bipartite quiver
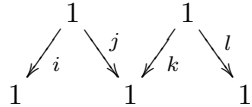
\[
\begin{xy}
\xymatrix@R2pt@C10pt{
&1\ar[lddd]^{i}\ar[rddd]^{j}&&1\ar[lddd]^{k}\ar[rddd]^{l}&\\\\\\1&&1&&1}
\end{xy}
\]
we always assumed $c(j)\neq c(k)$. But if we consider the quiver
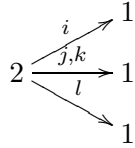
\[
\begin{xy}
\xymatrix@R10pt@C30pt{
&1\\2\ar[ru]^{i}\ar[rd]^l\ar[r]^{j,k}&1\\&1}
\end{xy}
\]
we could in a sense understand this quiver as the case $c(j)=c(k)$. But this raises another problem: we get additional conditions for $c(i)$ and $c(l)$ and moreover different symmetries. For instance, in the first case we have the symmetries of $S_2$. But in the second one we have the symmetries of $S_3$.\\

Another question is how to count or get all localization data (at least all of type one). Unfortunately, by use of the glueing method we do not get all localization data of type one. If it were possible to get all data of this type and if it could be shown that the other data come in a way from quivers of type one, one could probably prove the continuity. This would suffice to prove the existence of the conjectured function.\\

Finally, we give an example for a quiver of type one, which cannot be constructed by use of the glueing method. Let $(d,e)=(7,10)=(2,3)+(5,7)=(2,3)+(2,3)+(3,4)$. Then we have $(d_s,e_s)=(2,3)$. Consider
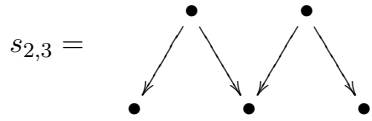
\[
\begin{xy}
\xymatrix@R2pt@C10pt{
&&\bullet\ar[lddd]\ar[rddd]&&\bullet\ar[lddd]\ar[rddd]&\\s_{2,3}=\\\\&\bullet&&\bullet&&\bullet}
\end{xy}
\]
and
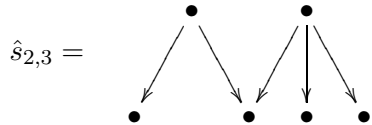
\[
\begin{xy}
\xymatrix@R2pt@C10pt{
&&\bullet\ar[lddd]\ar[rddd]&&\bullet\ar[lddd]\ar[ddd]\ar[rddd]&\\\hat{s}_{2,3}=\\\\&\bullet&&\bullet&\bullet&\bullet}
\end{xy}
\]
where the dots represent vertices of dimension one. We get the data of dimension type $(3,4)$ in the same way. But we do not get the following localization data of dimension type $(9,13)=(2,3)+(7,10)$ by sticking together the above ones:

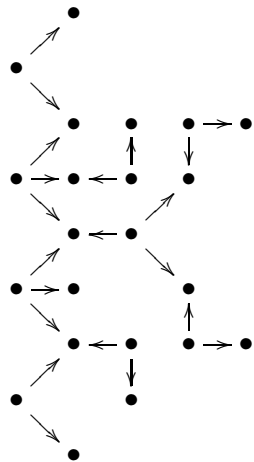
\[
\begin{xy}
\xymatrix@R10pt@C10pt{
&\bullet\\\bullet\ar[ru]\ar[rd]&\\&\bullet&\bullet&\bullet\ar[d]\ar[r]&\bullet\\\bullet\ar[ru]\ar[r]\ar[rd]&\bullet&\bullet\ar[l]\ar[u]&\bullet\\&\bullet&\bullet\ar[l]\ar[ru]\ar[rd]\\
\bullet\ar[ru]\ar[r]\ar[rd]&\bullet&&\bullet\\
&\bullet&\bullet\ar[l]\ar[d]&\bullet\ar[u]\ar[r]&\bullet\\\bullet\ar[ru]\ar[rd]&&\bullet\\&\bullet}
\end{xy}
\]
\section*{Acknowledgements}\noindent I would like to thank Markus Reineke for his support and for very helpful discussions.

\end{document}